\documentclass[11pt]{amsart}

\usepackage{amssymb, amsmath, amscd, amsthm,
}
\usepackage[all]{xy}          
\xyoption{dvips}              

\newcommand{\R}{{\mathbb R}}

\newcommand{\Q}{{\mathbb Q}}

\newcommand{\Field}{{\mathbb K}}

\newcommand{\BbbQ}{{\mathbb Q}}
\newcommand{\BbbR}{{\mathbb R}}
\newcommand{\BbbN}{{\mathbb N}}
\newcommand{\StrictOO}{\mathcal{O}^{s}}
\newcommand{\HH}{\operatorname{H}} 

\newcommand{\E}{\Emb}
\newcommand{\Emb}{\operatorname{Emb}}
\newcommand{\Ebar}{\overline{\Emb}}
\newcommand{\LEmb}{\operatorname{AffEmb}}

\newcommand{\calJ}{{\mathcal J}}
\newcommand{\LImm}{\operatorname{AffImm}}
\newcommand{\LEbar}{\overline{\LEmb}}
\newcommand{\Imm}{\text{Imm}}
\newcommand{\Spectra}{\operatorname{Spectra}}

\newcommand{\map}{\operatorname{Map}}

\newcommand{\holim}{\operatorname{holim}}

\newcommand{\hofiber}{\operatorname{hofiber}}
\newcommand{\HQ}{\mathbf{H}{\mathbb Q}}
\newcommand{\Top}{\operatorname{Top}}
\newcommand{\CC}{\operatorname{Ch}}
\newcommand{\Ch}{\operatorname{Ch}}

\newcommand{\calA}{\mathcal{A}}
\newcommand{\calB}{\mathcal{B}}
\newcommand{\calC}{\mathcal{C}}
\newcommand{\calD}{\mathcal{D}}
\newcommand{\calE}{\mathcal{E}}
\newcommand{\calO}{\mathcal{O}}
\newcommand{\Os}{\mathcal{O}^s}
\newcommand{\tOs}{\widetilde{\mathcal{O}^s}}
\newcommand{\tOsk}{\widetilde{\mathcal{O}^s_k}}

\newcommand{\calR}{\mathcal{R}}
\newcommand{\Balls}{\operatorname{B}}

\newcommand{\Config}{\mbox{C}}
\renewcommand{\O}{\mathcal{O}}
\newcommand{\OO}{\operatorname{O}}

\newcommand{\iso}{\operatorname{Iso}}

\newcommand{\Sing}{\operatorname{C}_*}

\newcommand{\Po}{\operatorname{Po}}
\newcommand{\Normalized}{\operatorname{N}}

\newcommand{\Vectors}{{\mathcal V}}
\newcommand{\sQmod}{{\Vectors}^{\Delta^{\operatorname{op}}}}
\newcommand{\CI}{{\mathcal C}I}
\newcommand{\ChI}{{\Ch}I}

\newcommand{\unit}{{\mathbf{1}}}
\newcommand{\rmods}[1]{{\mathcal Mod}\!\!-\!\!{#1}}

\newcommand{\quism}{\stackrel{\simeq}{\to}}
\newcommand{\Set}{\operatorname{Set}}
\newcommand{\ChImodule}{{$\ChI$-module}}
\newcommand{\ed}{\operatorname{ED}}
\theoremstyle{plain}
\newtheorem{thm}{Theorem}[section]
\newtheorem{prop}[thm]{Proposition}
\newtheorem{lemma}[thm]{Lemma}
\newtheorem{cor}[thm]{Corollary}
\newtheorem{conj}[thm]{Conjecture}

\theoremstyle{definition}
\newtheorem{definition}[thm]{Definition}
\newtheorem{eg}{Example}

\theoremstyle{remark}
\newtheorem{rem}[thm]{Remark}

\newcommand{\refT}[1]{Theorem~\ref{T:#1}}
\newcommand{\refC}[1]{Corollary~\ref{C:#1}}

\newcommand{\refP}[1]{Proposition~\ref{P:#1}}
\newcommand{\refD}[1]{Definition~\ref{D:#1}}
\newcommand{\refL}[1]{Lemma~\ref{L:#1}}

\newcommand{\refS}[1]{Section~\ref{S:#1}}

\begin{document}


\title[Calculus, formality, and embedding spaces]{Calculus of functors, operad formality, and rational homology of embedding spaces}


\author{Gregory Arone}
\address{Department of Mathematics, University of Virginia,
Charlottesville, VA} \email{zga2m@virginia.edu}
\author{Pascal Lambrechts}
\address{Institut Math\'{e}matique, Universit\'e de Louvain, 2 Chemin du Cyclotron, B-1348 Louvain-la-Neuve, Belgium}
\email{lambrechts@math.ucl.ac.be}
\author{Ismar Voli\'c}
\address{Department of Mathematics, Wellesley College,
Wellesley, MA} \email{ivolic@wellesley.edu}
\thanks{The first and third authors were supported by the National Science Foundation, 
grants DMS 0605073 and DMS 0504390 respectively. The second author is chercheur qualifi\'e au
F.N.R.S. and he gratefully acknowledges support by the Institut Mittag-Leffler (Djursholm, Sweden).}


\begin{abstract}

Let $M$ be a smooth manifold and $V$ a Euclidean space. Let $\Ebar(M,V)$ be the homotopy fiber of the map $\Emb(M,V)\longrightarrow \Imm(M,V)$. This paper is about the
rational homology of $\Ebar(M,V)$. We study it by applying 
embedding calculus and orthogonal calculus to the bi-functor $(M,V)\mapsto
\HQ\wedge\Ebar(M,V)_+$.  Our main theorem states that if $\dim V\geq 2\ed(M)+1$ (where $\ed(M)$ is the embedding dimension of $M$), the Taylor tower in the sense of orthogonal calculus (henceforward called ``the orthogonal tower'') of this functor  splits as a product of its layers.  Equivalently,
the rational homology spectral sequence associated with the tower
collapses at $E^1$. In the case of knot embeddings, this spectral
sequence coincides with the Vassiliev spectral sequence. The main
ingredients in the proof are embedding calculus and
Kontsevich's theorem on the formality of the little balls operad.

We write explicit formulas for the layers in the orthogonal tower
of the functor $$\HQ\wedge\Ebar(M,V)_+.$$ The formulas show, in
particular, that the (rational) homotopy type of the layers of the
orthogonal tower is determined by the (rational) homotopy type of
$M$. This, together with our rational splitting theorem, implies
that, under the above assumption on codimension,  rational homology equivalences 
of manifolds induce isomorphisms between the rational homology groups of 
$\Ebar(-,V)$.

\end{abstract}

\maketitle

\tableofcontents


\section{Introduction}\label{S:Intro}


Let $M$ be a smooth manifold of
dimension $m$. $M$ may be non-compact, 
but we always assume that $M$ is the interior of a compact manifold with boundary. Let $V$ be a Euclidean space. Let 
$\E(M,V)$ be the space of smooth embeddings of $M$ into $V$. For 
technical reasons, rather than study $\E(M,V)$ directly, we will focus on the space
$$\Ebar(M,V):=\hofiber\left(\E(M,V)\longrightarrow \Imm(M,V)\right),$$ where
$\Imm(M,V)$ denotes the space of immersions of $M$ into $V$. Note that the
definition requires that we fix an embedding (or at minimum an
immersion) $\alpha:M\hookrightarrow V$, to act as a basepoint.
Most of the time we will work with the suspension spectrum
$\Sigma^\infty \Ebar(M,V)_+$, and our results are really about the
rationalization of this spectrum, $\Sigma^\infty_{\Q} \Ebar(M,V)_+
\simeq \HQ\wedge \Ebar(M,V)_+$. In other words, our results are
about the rational homology of $\Ebar(M,V)$.

Our framework is provided by the Goodwillie-Weiss calculus of functors. One of the main features of calculus of functors is that it associates to a functor a tower of fibrations, analogous to the Taylor series  of a function. 
The functor $\Emb(M,V)$ is a functor of two variables, and accordingly 
one may do ``Taylor expansion'' in at least two ways:  In either the variable $M$ or the variable $V$
(or both). Since the two variables of $\Emb(M,V)$ are of rather different 
nature (for example, one is contravariant and the other one is covariant), there
are two versions of calculus needed for dealing with them -- {\it embedding
calculus} (for the variable $M$) and {\it orthogonal calculus} (for the variable $V$).

Embedding calculus
\cite{WeissEmb,GW} is designed for studying contravariant isotopy
functors on manifolds, such as $F(M)=\E(M,V)$. 
To such a functor $F$, embedding calculus associates a tower of
fibrations under $F$
\begin{equation}\label{E:EmbeddingTower}
F(-)\longrightarrow \big( T_\infty
F(-)\longrightarrow\cdots\longrightarrow T_{k}F(-) \longrightarrow
T_{k-1}F(-)\longrightarrow\cdots\longrightarrow T_1 F(-)\big).
\end{equation}
Here
$$
T_kF(U):=\underset{\{U'\in\mathcal{O}_k(M)\mid U'\subset
U\}}{\mbox{holim}}\, F(U'),
$$
where $\mathcal{O}_k(M)$ is the category of open subsets of $M$ that
are homeomorphic to the disjoint union of at most $k$ open balls.

$T_\infty F$ is defined to be the homotopy inverse limit of
$T_kF$. When circumstances are favorable, the natural map $F(M)\to
T_\infty F(M)$ is a homotopy equivalence, and then one says that
the embedding tower converges. There is a deep and important
convergence result, due to Goodwillie and Klein (unpublished, see
\cite{GK}), for the functor $F(M)=\E(M,N)$, where $N$ is a fixed
manifold.  We will state it now, it being an important fact in the
background, but we will not really use it in this paper.
\begin{thm}[Goodwillie-Klein, \cite{GK}]\label{T: Hard Convergence}
The Taylor tower (as defined above) of the embedding functor
$\E(M,N)$ (or $\Ebar(M,N)$) converges if $\dim(N)-\dim(M)\ge 3$.
\end{thm}
We will only need a much weaker convergence result, whose proof is
accordingly easier. The ``weak convergence theorem'' says that the
above Taylor tower converges if $2\dim(M)+2<\dim(N)$ and a proof
can be found in the remark after Corollary 4.2.4 in
\cite{GKW-survey}.
 The weak convergence result also holds for $\HQ\wedge\Ebar(M,N)_+$ by the main
result of \cite{WeissConvergence}.

Let us have a closer look at the functor $U\mapsto\HQ\wedge
\Ebar(U,V)_+$. If $U$ is homeomorphic
to a disjoint union of finitely many open balls, say $U\cong
k_U\times D^m$, then $\Ebar(U,V)$ is homotopy equivalent to the
configuration space $\Config(k_U,V)$ of $k_U$-tuples of distinct
points in $V$ or, equivalently, the space of $k_U$-tuples of
disjoint balls inside the unit ball of $V$, which we denote $\Balls(k_U,V)$.  Abusing
notation slightly, we can write that
\begin{equation}\label{TkEbar}
T_k\HQ\wedge\Ebar(M,V)_+:=\underset{U\in\mathcal{O}_k(M)}\holim
\HQ\wedge \Ebar(U,V)_+
\simeq\underset{U\in\mathcal{O}_k(M)}{\mbox{holim}}\, \HQ\wedge
\Balls(k_U,V)_+
\end{equation}
The right hand side in the above formula is not really
well-defined, because $\Balls(k_U,V)$ is not a functor on
$\mathcal{O}_k(M)$, but it gives the right idea. The formula tells
us that under favorable circumstances (e.g., if
$2\dim(M)+2<\dim(V)$), the spectrum $\HQ\wedge \Ebar(M,V)_+$ can
be written as a homotopy inverse limit of spectra of the form
$\HQ\wedge \Balls(k_U,V)_+$. It is obvious that the maps in the
diagram are closely related to the structure map in the little
balls operad. Therefore, information about the rational homotopy
type of the little balls operad may yield information about the
homotopy type of spaces of embeddings. The key fact about the
little balls operad that we want to use is the theorem of
Kontsevich (\cite[Theorem 2 in Section 3.2]{K_Formality}), asserting that this operad is  formal. Let $\Balls(\bullet, V)=\{\Balls(n, V)\}_{n\geq 0}$ be our notation for the operads of little balls inside the unit ball of $V$.

\begin{thm}[Kontsevich, \cite{K_Formality}] \label{T:KontsevichFormality}
The little balls operad is 
formal over the reals.  Explicitly, there is a chain of quasi-isomorphisms of
operads of chain complexes connecting the operads $\Sing(\Balls
(\bullet,V))\otimes \R$ and $\HH_*(\Balls (\bullet,V);\R).$
\end{thm}

The formality theorem was announced by Kontsevich in \cite{K_Formality}, and an outline of a proof was given there. However, not all the steps of the proof are given in \cite{K_Formality} 
in as much detail as some readers might perhaps wish. Because of this, the second and the third author
decided to write another paper \cite{LV_Formality}, whose primary purpose is to spell out the details of the proof of the formality theorem, following Kontsevich's outline. The paper \cite{LV_Formality} also has a second purpose, which is to prove a slight strengthening of the formality theorem, which we
call ``a relative version'' of the formality theorem (Theorem
\ref{T:relativeKontsevich} in the paper). We will give a sketch of the proof of the relative version in \refS{formalityballs}. Using the relative version of formality, together with some abstract homotopy
theory, we deduce our first theorem (see
\refT{formalityembeddingdiagram} for a precise statement).
\begin{thm}\label{T:formality}
Suppose that the basepoint embedding $\alpha:M\hookrightarrow V$
factors through a vector subspace $W\subset V$ such that
$\dim(V)\ge 2\dim(W)+1$. Then the contravariant functor from $\tOsk(M)$ to chain complexes
$$U\mapsto \Sing(\Ebar(U,V))\otimes \R$$
 is formal. This means that there is a chain of weak equivalences, natural in $U$
$$\Sing(\Ebar(U,V))\otimes \R \simeq \HH_*(\Ebar(U,V);\R)$$
\end{thm}
Here, $\tOsk(M)$ is a suitable variation of the category $\mathcal{O}_k(M)$ (where $k$ can be arbitrarily large). We will now give the rough idea of the proof. The category $\tOsk(M)$ is a subcategory of $\tOsk(W)$, so it is enough to prove that $\Sing(\Ebar(U,V))\otimes \R$ is formal as a functor on $\tOsk(W)$ (this is Theorem~\ref{T:formalityembeddingdiagram} in the paper). The category $\tOsk(W)$, is a category of balls in $W$, and so it is related to the little balls operad $\Balls(\bullet, W)$. Therefore, the category of (contravariant) functors from $\tOsk(W)$ to (real) chain complexes is closely related to the category of (right) modules over the chains on little balls operad $\Sing(\Balls(\bullet,W))\otimes\R$. The space $\Ebar(U,V)$, where $U$ is the union of $n$ balls, is equivalent to the $n$-th space in the 
$V$-balls operad, $\Balls(n,V)$. We will show that the formality of the functor $\Sing(\Ebar(U,V))\otimes \R$ follows from the formality of the operad $\Sing(\Balls(\bullet,V))\otimes \R$ {\em as a module over the operad $\Sing(\Balls(\bullet,W))\otimes \R$}. The last formality statement follows from the relative formality theorem. 
\begin{rem}\label{R: Caveat}
In several places, our argument relies on the following simple observation: the formality of a discrete (i.e., unenriched) diagram of chain complexes is equivalent to the splitting up to homotopy of the  diagram as a direct sum of diagrams concentrated in a single homological degree (Proposition~\ref {P:PostnikovInterpretation}). This is convenient, because this kind of splitting is a homotopy invariant property, that is preserved by various Quillen equivalences between diagram categories that we need to consider, while the property of being formal can not, in general, be transferred across a Quillen equivalence. 
However, for enriched diagrams, this observation is not true. We can see it in the example of modules over operads (which are a special case of enriched functors). The homology the little balls operad is formal, but, plainly, it does not split, as a module over itself, into a direct sum of modules concentrated in a single homological dimension. However, if we consider the operad $\HH_*(\Balls(\bullet,V))$ as a module over the operad $\HH_*(\Balls(\bullet, W))$, where $W$ is a proper subspace of $V$, then it does split as a direct sum, for the silly reason that in this case the action of $\HH_*(\Balls(\bullet, W))$ factors through an action of $\HH_0(\Balls(\bullet, W))$, and so it is essentially a discrete action. This follows from the elementary fact that a proper linear inclusion of vector spaces $W\hookrightarrow V$ induces a null-homotopic map of spaces of little balls $\Balls(n, W)\longrightarrow \Balls(n, V)$. Using this observation, together with the relative formality theorem, we are able to conclude the highly non-obvious statement that if $2\dim(W)<\dim(V)$ then the operad $\Sing(\Balls(\bullet,V))\otimes \R$ splits, as a module over $\Sing(\Balls(\bullet,W))\otimes \R$, into a direct sum of modules that are homologically concentrated in a single dimension. This splitting is then used to prove an analogous splitting (and therefore the formality) of the functor $\Sing(\Ebar(U,V))\otimes \R$. 
It is for this reason that we have to assume that $M$ lies in a proper linear subspace of $V$, and to utilize a relative formality theorem.
\end{rem}

A formality theorem similar to \refT{formality} was used in
\cite{LV} for showing the collapse (at $E^2$) of a certain
spectral sequence associated to the {\em embedding} tower for
spaces of \emph{knot embeddings}.  However, to obtain a collapsing result
for a spectral sequence for more general embedding spaces, we need,
curiously enough, to turn to Weiss' {\em orthogonal} calculus
(the standard reference is \cite{WeissOrth}, and a brief overview can be found in Section \ref{S:ReviewCalculi}).  This is a calculus of covariant
functors from the category of vector spaces and linear isometric
inclusions to topological spaces (or spectra). To such a functor
$G$, orthogonal calculus associates a tower of fibrations of
functors $P_nG(V)$, where $P_nG$ is the $n$-th Taylor polynomial
of G in the orthogonal sense. Let $D_nG(V)$ denote the $n$-th
homogeneous layer in the orthogonal Taylor tower, namely the fiber
of the map $P_nG(V)\to P_{n-1}G(V)$.

The functor that we care about is, of course, $G(V)=\HQ\wedge
\Ebar(M,V)_+$ where $M$ is fixed. We will use the notation
$P_n\HQ\wedge \Ebar(M,V)_+$ and $D_n\HQ\wedge \Ebar(M,V)_+$ to
denote its Taylor approximations and homogeneous layers in the
sense of orthogonal calculus. It turns out that \refT{formality}
implies that, under the same condition on the codimension, the
orthogonal tower of $\HQ\wedge \Ebar(M,V)_+$ splits as a product
of its layers. The following is our main theorem (\refT{TowerSplits} in the paper).
 \begin{thm}\label{T:OrthogonalSplits}
Under the assumptions of \refT{formality}, there is a homotopy
equivalence, natural with respect to embeddings in the
$M$-variable (note that we do not claim that the splitting is natural in $V$)
$$P_n\HQ\wedge\Ebar(M,V)_+\simeq \prod_{i=0}^n D_i\HQ\wedge\Ebar(M,V)_+.$$
\end{thm}
The following corollary is just a reformulation of the theorem.
\begin{cor} Under the assumptions of \refT{formality} and \refT{OrthogonalSplits},
the spectral sequence for $\HH_*(\Ebar(M,V);\Q)$ that arises from
the Taylor tower (in the sense of orthogonal calculus) of
$\HQ\wedge\Ebar(M,V)_+$ collapses at $E^1$.
\end{cor}
Here is a sketch of the proof of \refT{OrthogonalSplits}.
Embedding calculus tells us, roughly speaking, that $\HQ\wedge \Ebar(M,V)_+$ can be
written as a homotopy limit of a diagram of spectra of the from
$\HQ\wedge \Balls(k,V)_+$. Since there is a Quillen equivalence between
the categories of rational spectra and rational chain complexes, we may pass to a
diagram of rational chain complexes of the form $\Sing(\Balls(k,V))\otimes \Q$, whose homotopy 
limit is $\Sing(\Ebar(M,V))\otimes \Q$.   
On the other hand, \refT{formality} tells us
that this diagram of chain complexes is formal when tensored with $\R$. One concludes that $\Sing( \Ebar(M,V))\otimes \R$ splits as the product of inverse limits of layers in the 
Postnikov towers of $\Sing(\Balls(k,V))\otimes \R$. It turns out that in our case tensoring with $\R$ commutes with taking the homotopy limit, and so it follows that there must be a similar splitting for $\Sing( \Ebar(M,V))\otimes \Q$ and therefore for 
$\HQ\wedge \Ebar(M,V)_+$. On the other hand, it
turns out that for functors of the form $\HQ\wedge
\Balls(k,V)_+$, the Postnikov tower coincides, up to regrading,
with the orthogonal tower, and
therefore $\HQ\wedge \Ebar(M,V)_+$ splits as the product of
inverse limits of layers in the orthogonal tower of rationally stabilized 
configuration spaces. But, taking the $n$-th layer in
the orthogonal tower is an operation that commutes
 (in our case) with homotopy
inverse limits (unlike the operation of taking the $n$-th layer of the Postnikov
tower), and therefore $\HQ\wedge \Ebar(M,V)_+$ splits as the
product of layers of its orthogonal tower.
\begin{rem} In the case of knot embeddings, the spectral sequence
associated with the orthogonal tower seems to coincide with the famous spectral
 sequence constructed by Vassiliev, since the latter also collapses, and the initial terms are isomorphic. We hope to come back to this in future papers.
\end{rem}
In \refS{DescriptionOfLayers} we write an explicit description of
$D_n\Sigma^\infty\Ebar(M,V)_+$, in terms of 
certain spaces of partitions (which can also be described as
spaces of rooted trees) attached to $M$. This section is an
announcement; detailed proofs will appear in \cite{Arone}. We do
note the following consequence of our description of the layers:  The 
homotopy groups of the layers depend only on the stable homotopy type of $M$
 and similarly the rational homotopy groups of the layers depend only on the rational homology type of $M$ (\refC{RHT-invariance-layers}). Combining this with
\refT{OrthogonalSplits}, we obtain the following theorem (\refT{RationalHomologyInvariant} in the paper).
\begin{thm}\label{T:MainCorollary} Under the assumptions of \refT{OrthogonalSplits},
the rational homology groups of the space $\Ebar(M,V)$ are
determined by the rational homology type of $M$. More precisely, suppose $M_1,M_2, V$ satisfy the assumptions of \refT{OrthogonalSplits}, and suppose there is a zig-zag of maps, each inducing an isomorphism in rational homology, connecting $M_1$ and $M_2$. Then there is an isomorphism
$$\HH_*(\Ebar(M_1,V);\BbbQ)\cong \HH_*(\Ebar(M_2,V);\BbbQ).$$
\end{thm}
In view of this result, one may wonder whether the rational {\em homotopy} groups, or may be even the rational homotopy \emph{type} (rather than just rational homology) of
$\Ebar(M,V)$ could be an invariant of the rational homotopy type
of $M$ (in high enough codimension). One could derive further hope from the fact
that the little balls operad is not only formal, but also coformal. We would like to propose the following conjecture
\begin{conj} Under the assumptions of \refT{formality},
the rational homotopy spectral sequence for $\pi_*(\Ebar(M,V))\otimes \Q$ that arises from
the Taylor tower (in the sense of orthogonal calculus) of
$\Ebar(M,V)$ collapses at $E^1$.
\end{conj}
A statement essentially equivalent to this conjecture is proved
in \cite{FourAuthors} in the special case of spaces of long knots in dimension
$\ge 4$.
In this case, the space of embeddings is an $H$-space (in fact,
a double loop space), and this gives one enough control over 
the homotopy type of the space to force the desired conclusion.

A general point that we are trying to make with this paper is this: while embedding
calculus is important, and is in some ways easier to understand
than orthogonal calculus, the Taylor tower in the sense of
orthogonal calculus is also interesting and is worthy of a further
study. We hope that \refS{DescriptionOfLayers} will convince the
reader that the layers of the orthogonal tower, while not exactly
simple, are interesting, and it may be possible to do calculations
with them. We intend to come back to this in the future.

\subsection{A section by section outline}
In \refS{SpacesSpectraChainComplexes} we review background
material and fix terminology on spaces, spectra and chain
complexes.  In \refS{Formality} we define the notion of 
formality of diagrams chain complexes. The main
result of this section is the following simple but useful observation: the stable
formality of a diagram can be interpreted as the splitting of its Postnikov tower.

Our next goal is to exploit  Kontsevich's formality of the little
balls operads and deduce some formality results of diagrams of
embedding spaces. In order to do that we first review, in
\refS{enriched}, enriched categories, their modules and the
associated homotopy theory. In Section \ref{S:operads} we review
classical operads and their modules and give an alternative
viewpoint on those in terms of enriched categories. This will be
useful for the study of the homotopy theory of modules over an
operad. We then digress in Section
\ref{S:formalityballs} to prove a relative version of Kontsevich's
formality of the little balls operads that we need for our applications. 
In \refS{FormalityEmbedding} we
deduce the formality of a certain diagram of real-valued chains on embedding spaces.

In Section \ref{S:ReviewCalculi} we digress again to give a review of embedding calculus and orthogonal calculus, and record some generalities on how these two brands
of calculus may interact. In Section \ref{S:RatEmb} we use the formality of a diagram of 
chains on embedding spaces established in \refS{FormalityEmbedding}
to show that the stages in the embedding tower of $\HQ\wedge \Ebar(M,V)_+$ split in a certain way, but \emph{not} as the product of the layers in the embedding tower. In Section \ref{S:FormalityOrthogonal} we reinterpret this splitting once again, to prove our main theorem:  Under a certain co-dimension hypotheses, the orthogonal tower of
$\HQ\wedge\Ebar(M,V)_+$ splits as the product of its layers. In Section \ref{S:DescriptionOfLayers} we sketch a description of the layers in the orthogonal tower, and deduce
that the rational homology of the space of embeddings (modulo immersions) of a manifold into a
high-dimensional vector space is determined by the rational
homology type of the manifold.


\subsection{Acknowledgments} The second author thanks Enrico Vitale for help with enriched categories.



\section{Spaces, spectra, and chain complexes}\label{S:SpacesSpectraChainComplexes}

\noindent Let us introduce the basic categories that we will work with.

\textbullet \ \ $\Top$ will stand for the category of compactly generated spaces
(we choose \emph{compactly generated} to make  it  a \emph{closed}
monoidal category, see Section \ref{S:enriched}). If $X$ is a
space we denote by $X_+$ the based space obtained by adjoining a
disjoint basepoint.

\textbullet \ \ $\Spectra$ will be the category of $(-1)$-connected spectra. We denote by
$\HQ$ the Eilenberg-MacLane spectrum such that $\pi_0(\HQ)=\BbbQ$.
A rational spectrum is a module spectrum over $\HQ$. For a space $X$, $\Sigma^\infty X_+$
stands for the suspension spectrum of $X$, and $\HQ\wedge X_+$ denotes
the stable rationalization of $X$. It is well-known that there is a \emph{rational} equivalence
$\HQ\wedge X_+\simeq\Sigma^\infty X_+$. 

\textbullet \ \ $\Vectors$ will denote the category of rational vector spaces (or $\Q$-vector spaces), and
$\sQmod$ the category of simplicial $\Q$-vector spaces.

\textbullet \ \ $\CC_{\Q}$ and $\CC_{\R}$ will denote the category of non-negatively graded rational and real chain complexes respectively. We will some times use $\CC$ to denote either one of these
two categories.

\vskip 6pt
\noindent
Most of the above categories have a {\em Quillen model structure}, 
which means that one can apply to them the techniques of homotopy
theory. A good introduction to closed model categories is
\cite{DS}, a good reference is \cite{Hir}. There are slight variations in the literature
as to the precise definition of model structure. We use the definition given in \cite{Hir}. In particular, we assume the existence of functorial fibrant and cofibrant replacements. The category for which we will 
use the model structure most heavily is the category of chain complexes.
Thus we remind the reader that the category of chain
complexes over a field has a model structure where weak equivalences are
quasi-isomorphisms, fibrations are chain maps that are surjective
in positive degrees, and cofibrations are (since all modules are projective) chain maps that are injective in all degrees \cite[Theorem 7.2]{DS}. We will also need the fact that the
category of rational spectra is a Quillen model category and is Quillen equivalent
to the category $\CC_{\Q}$. For a proof of this (in fact, of a more general statement, involving the category of module spectra over a general Eilenberg - Mac Lane commutative ring-spectrum) see, for example, \cite{Schwede}.

We now define some basic functors between the various
categories in which we want to do  homotopy theory. 

\subsubsection{Homology}\label{S:Homology}

We think of homology as a functor from chain complexes to chain
complexes. Thus if $C$ is a chain complex, then $\HH_*(C)$ is the chain complex whose chain groups are the homology groups of $C$, and whose differentials are zero. 
Moreover, we define $\HH_n(C)$
to be the chain complex having the $n$-th homology group of $C$ in degree $n$ and zero in all other degrees. Thus, $\HH_n$ is a functor from $\CC$ to $\CC$ as well. Notice that there are obvious isomorphisms of functors
$$\HH_*\cong \bigoplus_{n=0}^\infty \HH_n \cong \prod_{n=0}^\infty \HH_n.$$

\subsubsection{The normalized chains functor}\label{S:Dold-KanCorrespondence}


 To get from
spaces to chain complexes, we will use the {\em normalized 
singular chains} functor $\Sing:\Top\to\CC$, defined as 
$$\Sing(X)=\Normalized(\Q[{\mathcal S}_\bullet(X)]).$$
Here ${\mathcal S}_\bullet(X)$ is the simplicial set of singular simplices of $X$,
$\Q[{\mathcal S}_\bullet(X)]$ is the simplicial $\Q$-vector
space generated by ${\mathcal S}_\bullet(X)$, and
$\Normalized:\sQmod\to \CC$ is the normalized chains functor as
defined for example in \cite[Chapter 8]{Weibel}.

\subsection{Postnikov sections}\label{S:PostnikovSections}


We will need to use Postnikov towers in the categories of chain complexes, and spectra. We now review  the construction of Postnikov towers in the category of chain complexes.
For an integer $n$ and a chain complex $(C,d)$, let 
$d(C_{n+1})$ be the $n$-dimensional boundaries in $C$. We define the
$n$th-Postnikov section of $C$, denoted 
$(\Po_n(C),d')$,  as follows
$$ \left(\Po_n(C)\right)_i=
  \begin{cases}
    C_i & \text{if $i\leq n$}, \\
    d(C_{n+1}) & \text{if $i=n+1$},\\
    0 & \text{if $i>n+1$},
  \end{cases}
$$
The differential $d'$ is defined to be $d$ in degrees $\le n$, and the obvious inclusion $d(C_{n+1})\hookrightarrow C_n$ in degree $n+1$. It is easy to see that $\Po_n$ defines  a
functor from $\CC$ to $\CC$. Moreover, $\HH_i(\Po_n(C))\cong
\HH_i(C)$ for $i\leq n$ and $\HH_i(\Po_n(C))=0$ for $i>n$.

For each $n$, there is a natural fibration (i.e., a degree-wise surjection)
$\pi_n\colon\Po_{n}(C)\twoheadrightarrow\Po_{n-1}(C)$ defined as follows: $\pi_n$ is the identity in all degrees except $n+1$ and $n$; in degree $n+1$ it is the zero homomorphism; and in degree $n$ it is the obvious surjective map $d:C_n\to d(C_n)$. Since $\pi_n$ is a fibration, $\ker(\pi_n)$ can serve as the model for its homotopy fiber. Clearly, $\ker(\pi_n)$ is a chain complex concentrated in dimensions $n$ and $n+1$. The homology of the kernel is concentrated in dimension $n$, and in this dimension it equals the homology of the original complex $C$.
A similar formula defines a natural map $\rho_n\colon
C\to\Po_n(C)$, and we have $\pi_n\rho_{n}=\rho_{n-1}$. Note that 
$\rho_n$, like $\pi_{n+1}$, is an isomorphism (on chain level) in degrees
$\leq n$.

\subsection{Diagrams}\label{S:Diagrams}

 Let  $\calA$ be a small
category and let $\calE$ be a category. An \emph{$\calA$-diagram
in $\calE$} is just a functor $F\colon\calA\to\calE$. In this
paper a diagram can be a functor which is either covariant or
contravariant. A \emph{morphism of $\calA$-diagrams} is a natural
transformation between two functors. Such a morphism is called a
\emph{weak equivalence} if it is a weak equivalence objectwise,
for a given notion of weak equivalence in the category $\calE$.  In practice, we will only
consider diagrams of spaces, chain 
complexes or spectra.

\subsection{Homotopy limits} We will make heavy use of homotopy limits of diagrams in $\Spectra$ and in $\CC$.
 Homotopy limits of diagrams in a general model category are treated in \cite{Hir}, Chapter 19. 
 Generally, when we take the homotopy limit of a diagram, we assume that all the objects in the diagram are fibrant and cofibrant 
 - this will ensure ``correct'' homotopical behavior in all cases.
Since most of our homotopy limits will be taken in the category of chain complexes over $\Q$ or $\R$, in which all objects are fibrant and cofibrant,
this is a moot point in many cases. The only other category in which we will take homotopy limits is the category of rational spectra, 
in which case we generally assume that we have taken fibrant-cofibrant replacement of all objects, whenever necessary. 

It follows from the results in \cite{Hir}, Section 19.4, that if $R$ and $L$ are the right and left adjoint 
in a Quillen equivalence, then both $R$ and $L$ commute with homotopy limits up to a zig-zag of natural weak equivalences. In particular, this enables us to shuttle back and forth between homotopy limits of diagrams of rational spectra and diagrams of rational chain complexes.

\section{Formality and homogeneous splitting of diagrams}\label{S:Formality}


The notion of formality was first introduced by Sullivan in the
context of rational homotopy theory
\cite{Sullivan-infinitesimal,DGMS}. Roughly speaking a chain
complex (possibly with additional structure) is called formal if
it is weakly equivalent to its homology. In this paper we will
only use the notion of formality of diagrams of chain complexes (over $\Q$ and over 
$\R$).

\begin{definition}
 Let ${\mathcal A}$ be a small category. An $\calA$-diagram of 
 chain complexes,  $F\colon{\mathcal A}\to\CC$,
 is \emph{formal} if there is a chain of weak equivalences
 $F\simeq \HH_*\circ F$.  
\end{definition}
Formality of chain complexes has a convenient interpretation as the splitting
of the Postnikov tower. 

\begin{definition}\label{D:discrete-splits}
Let $\calA$ be a small category. We say that an $\calA$-diagram of
chain complexes, $F\colon\calA\to\Ch$, \emph{splits homogeneously}
if there exist $\calA$-diagrams $\{F_n\}_{n\in\BbbN}$ 
of chain complexes such that $F\simeq \oplus_n
F_n$ and $\HH_*(F_n)=\HH_n(F_n)$ (i.e., $F_n$ is homologically concentrated in
degree $n$).
\end{definition}

\begin{prop}\label{P:PostnikovInterpretation}
Let $\mathcal A$ be a small category. An $\calA$-diagram of chain
complexes is formal if and only if it splits homogeneously.
\end{prop}

\begin{proof}
Let $F$ be an $\calA$-diagram of chain complexes.

In one direction, if $F$ is formal then $F\simeq \HH_*(F)$. Since
$\HH_*=\oplus_{n\in\BbbN}\HH_n$, we get the homogeneous splitting
$F\simeq \oplus_n \HH_n(F)$.

In the other direction, suppose that $F\simeq
\oplus_{n\in\BbbN}F_n$ with $\HH_*(F_n)=\HH_n(F_n)=\HH_n(F)$. Recall the
definition of Postnikov sections of chain complexes from Section
\ref{S:SpacesSpectraChainComplexes}. Then
$$
\ker\left(\Po_n(F_n) \stackrel{\pi_n}\twoheadrightarrow \Po_{n-1}(F_n)\right)
$$
is concentrated in degrees $n$ and $n+1$ and its homology is
exactly $\HH_n(F)$. Thus we have a chain of quasi-isomorphisms
$$F_n \stackrel{\simeq}{\longrightarrow} \Po_n(F_n) \stackrel{\simeq}{\longleftarrow} \ker\left(\Po_n(F_n) \to \Po_{n-1}(F_n)\right)
\stackrel{\simeq}{\longrightarrow} \HH_n\left(\ker(\Po_n(F_n) \to \Po_{n-1}(
F_n))\right)\cong \HH_n(F),$$ and so $F\simeq \oplus_n
\HH_n(F)=\HH_*(F)$.
\end{proof}
\begin{rem}\label{R:UniquenessEMdiagrams}
Note that in the above we proved the following (elementary) statement:  
Suppose $F$ and $G$ are two $\calA$-diagrams
of chain complexes such that both $F$ and $G$ are homologically concentrated in degree $n$ and such that there is an isomorphism of diagrams $\HH_n(F)\cong\HH_n(G)$. Then there is a chain of weak equivalences, $F\simeq G$. Using the Quillen equivalence between rational spectra and rational chain complexes, 
one can prove the analogous statement for diagrams of rational Eilenberg-Mac Lane spectra:  If $F$ and $G$ are two $\calA$-diagrams of
rational Eilenberg-Mac Lane spectra concentrated in degree $n$, and if there is an isomorphism of diagrams $\pi_n(F)\cong\pi_n(G)$
 then there is a chain of weak equivalences $F\simeq G$.
 Actually this is true for diagrams of (non rational) Eilenberg-MacLane spectra concentrated in the same degree but we will not need that fact.
\end{rem}

\begin{rem}\label{R:FormalLimitSplits}
Let $F$ be a diagram with values in $\CC$. There is a tower of fibrations converging to $\holim F$ whose $n$-th stage is $\holim \Po_nF$. We call it \emph{the lim-Postnikov tower}. Of course, this tower does not usually coincide with the Postnikov tower of $\holim F$. Since $\HH_*\cong \prod_{n=0}^\infty \HH_n$, and homotopy limits commute with products, it follows immediately that if $F$ is a formal diagram then the lim-Postnikov tower of $\holim  F$ splits as a product, namely $$\holim  F\simeq \prod_{n=0}^\infty \holim \HH_n\circ F$$ 
\end{rem}

The proof of the following is also straightforward.
\begin{lemma}\label{L:formalprecomposite}
Let $\lambda\colon\calA\to\calA'$ be a  functor between small
 categories and let $F$ be an $\calA'$-diagram of chain
complexes. If the
$\calA'$-diagram $F$ is formal then so is the
$\calA$-diagram $\lambda^*(F):=F\circ\lambda$.
\end{lemma}

\section{Enriched categories and their modules}\label{S:enriched}

We now briefly recall some definitions
and facts about symmetric monoidal categories, enriched categories, Quillen
model structures, etc. The standard reference for symmetric
monoidal categories and enriched categories is \cite[Chapter
6]{CategoryTheory}. We will also need some results of
Schwede and Shipley on the homotopy theory of enriched categories
developed in \cite{SS}, especially Section 6, which is where we also borrow some of our notation and terminology from.

\subsection{Monoidal model categories and enriched categories}\label{S:MonoidalModelandEnrichedCategories}


A \emph{closed symmetric monoidal category} is a triple
$(\calC,\otimes,\unit)$ such that $\otimes$ and $\unit$ endows the
category $\calC$ with a symmetric monoidal structure, and such that,
for each object $Y$, the endofunctor $-\otimes
Y\colon\calC\to\calC\,,\,X\mapsto X\otimes Y$ admits a right
adjoint denoted by $\calC(Y,-)\colon Z\mapsto\calC(Y,Z)$. It is 
customary to think of $\calC(Y,Z)$ as an ``internal mapping object''.
Throughout this section, $\calC$ stands for a closed symmetric monoidal
category.

A \emph{monoidal model category} is a closed symmetric monoidal
category  equipped with a compatible Quillen model structure (see
\cite[Definition 3.1]{SS} for a precise definition).

The only examples of monoidal model categories that we will
consider in this paper are
\begin{enumerate}
\item The category $(\Top,\times,*)$ of compactly generated topological spaces with cartesian product;
\item The category $(\Ch,\otimes,\Field)$ of non-negatively graded chain complexes over $\Field$ (where $\Field$ is $\Q$ or $\R$), with tensor product.
\end{enumerate}
The internal $\hom$ functor in the category $\Ch$ is defined as follows. Let $Y_*,Z_*$ be chain complexes. Then $\Ch(Y_*,Z_*)$ is the chain complex that in positive degrees $p>0$ is defined by 
$$\Ch_p(Y_*,Z_*)=\prod_{n=0}^\infty \hom(Y_n,Z_{n+p})$$while in degree zero, we have 
$$\Ch_0(Y_*,Z_*)=\{\mbox{chain homomorphisms from } Y_* \mbox{ to } Z_*\}.$$
The differential in $\Ch(Y_*,Z_*)$ is determined by the formula $D(\{f_n\})=\{d_Zf_n-(-1)^p f_{n-1}d_Y\}$, for $f_n\in \hom(Y_n,Z_{n+p})$.


\subsection{Enriched categories}\label{S:EnrichedCategories}

A category $\O$ \emph{enriched over $\calC$}, or
a \emph{$\calC$-category}, consists of a class $I$
(representing the \emph{objects} of $\calO$), and, for any objects
$i,j,k\in I$, a $\calC$-object $\O(i,j)$ (representing the
morphisms from $i$ to $j$ in $\calO$) and $\calC$-morphisms
$$\O(i,j)\otimes \O(j,k)\longrightarrow \O(i,k),\mbox{ and } 
\unit\to \O(i,i)
$$
(representing the composition of morphisms in $\calO$ and the identity
morphism on $i$). These structure morphisms are required to be associative and unital in 
the evident
sense.  Notice that a closed symmetric
monoidal cateory $\calC$ is enriched over itself since $\calC(Y,Z)$ is an object of $\calC$. Following \cite{SS}, we
use the term \emph{$\calC I$-category} to signify a category enriched over $\calC$, whose set of objects is $I$.

Let $\calO$ be a $\CI$-category and $\mathcal R$ be a category
enriched over $\mathcal C$. A (covariant) \emph{functor enriched
over $\calC$}, or \emph{$\calC$-functor} from $\O$ to $\calR$,
$$M\colon\calO\longrightarrow\calR,$$
consists of an $\calR$-object $M(i)$ for every $i\in I$, and of
morphisms in $\mathcal C$
$$M(i,j)\colon\O(i,j) \longrightarrow {\mathcal R}(M(i),M(j)),$$
for every $i,j\in I$, that are associative and unital. There is an
analogous notion of a \emph{contravariant $\calC$-functor}.

A \emph{natural transformation enriched over $\calC$}, $\Phi\colon
M\to M'$, between two $\calC$-functors $M,M'\colon\calO\to\calR$
consists of $\calC$-morphisms
$$
\Phi_i\colon\unit\longrightarrow\calR(M(i),M'(i))
$$
for every object $i$ of $\calO$, that satisfy the obvious
commutativity conditions for a natural transformation (see
\cite[6.2.4]{CategoryTheory}). Notice that if $\calR=\calC$ then  
a morphism $\Phi_i\colon\unit\to\calC(M(i),M'(i))$ is the same as the adjoint morphism 
$\Phi(i)\colon M(i)\to M'(i)$ in $\calC$.

For fixed $\calC$ and $I$, we consider the collection of $\CI$-categories as a category in its own right. A {morphism} of $\CI$-categories is an enriched functor that
is the \emph{identity} on the set of objects.

Suppose now that $\mathcal C$ is a monoidal model  category. In
particular, $\mathcal C$ is equipped with a notion of weak
equivalence. Then we say that a morphism $\Psi:\O\to \calR$ of
$\CI$-categories is a \emph{weak equivalence} if it is a weak
equivalence pointwise, i.e., if the map $\O(i,j)\to \calR(i,j)$ is
a weak equivalence in  $\mathcal C$ for all $i,j\in I$.


\subsection{Homotopy theory of right modules over enriched categories}\label{S:HomotopyTheoryOfModules}

For a $\CI$-category $\O$, a {\em (right) $\O$-module} is a
contravariant $\calC$-functor from $\O$ to $\mathcal C$.
Explicitly an $\O$-module $M$ consists of objects $M(i)$ in
$\calC$ for $i\in I$ and (since $\calC$ is a \emph{closed}
monoidal category and since it is enriched over itself) of $\calC$-morphisms
$$M(j)\otimes \O(i,j)\longrightarrow M(i)$$
which are associative and unital. A \emph{morphism} of $\O$-modules, 
$\Phi\colon M\to M'$, is an enriched natural transformation,
i.e., a collection of $\calC$-morphism $\Phi(i)\colon
M(i)\to M'(i)$ satisfying the usual naturality requirements. 
Such a morphism of $\O$-module is a \emph{weak
equivalence} if each $\Phi(i)$ is a weak equivalence in $\calC$.
We denote by $\rmods{\O}$ the category of right $\O$-modules and
natural transformations.

Let $\Psi:\O\to \calR$ be a morphism of $\CI$-categories. Clearly,
$\Psi$ induces a \emph{restriction of scalars} functor on module
categories
\begin{align*}
\Psi^*:\rmods{\calR} & \longrightarrow \rmods{\O}  \\ 
                  M  & \longmapsto M\circ \Psi.
\end{align*}
As explained in \cite[page 323]{SS}, the functor $\Psi^*$ has a
left adjoint functor $\Psi_*$, also denoted $-\otimes_{\O}
{\calR}$ (one can think of $\Psi_*$ as the left Kan extension).
Schwede and Shipley \cite[Theorem 6.1]{SS} prove that under some
technical hypotheses on $\mathcal C$, the category $\rmods{\O}$
has a Quillen model structure, and moreover, if $\Psi$ is a weak
equivalence of $\CI$-categories, then the pair $(\Psi^*,\Psi_*)$
induces a Quillen equivalence of module categories.

We will need this result in the case ${\mathcal C}=\CC$. In
keeping with our notation, we use $\CC\! I$-categories to denote
categories enriched over chain complexes, with object set $I$.
Note that the category of modules over a $\CC\! I$-categories
admits coproducts (i.e. direct sums).

\begin{thm}[Schwede-Shipley, \cite{SS}] \label{T:EquivalenceOfModuleCategories}\ 

(1) Let $\O$ be a $\CC \!I$-category. Then $\rmods{\O}$ has a cofibrantly generated Quillen model structure,
with fibrations and weak equivalences defined objectwise.

(2) Let $\Psi:\O\to \calR$ be a weak equivalence of $\CC\!
I$-categories. Then $(\Psi^*,\Psi_*)$ induce a Quillen equivalence
of the associated module categories.
\end{thm}
\begin{proof}
General conditions on $\mathcal C$ that guarantee the result are
given in \cite[Theorem 6.1]{SS}. It is straightforward to check
that the conditions are satisfied by the category of chain
complexes (the authors of \cite{SS} verify them for various
categories of spectra, and the verification for chain complexes is
strictly easier).
\end{proof}

Let $\O$ and $\calR$ be $\CI$-categories and let $M$ and $N$ be
right modules over $\O$ and $\calR$ respectively. A morphism of
pairs $(\O,M)\to (\calR,N)$ consists of a morphism of $\calC
I$-categories $\Psi\colon\O\to \calR$ and a morphism of
$\O$-modules $\Phi\colon M\to \Psi^*(N)$. The corresponding
category of pairs $(\O,M)$ is called the \emph{$\calC I$-module  category}.

A morphism $(\Psi,\Phi)$ in $\calC I$-module  is called a \emph{weak
equivalence} if both $\Psi$ and $\Phi$ are weak equivalences. Two
objects of $\calC I$-module  are called \emph{weakly equivalent} if
they are linked by a chain of weak equivalences, pointing in
either direction.

In our study of the formality of the little balls operad, we will
consider certain splittings of $\calO$-modules into direct sums. The
following homotopy invariance property of such a splitting will be important.
\begin{prop}\label{P:SplittingPropagates}
Let $(\O,M)$ and  $(\calO',M')$ be weakly equivalent \ChImodule s.
If $M$ is weakly equivalent as an $\calO$-module to a direct sum
$\oplus M_n$, then $M'$ is weakly equivalent as an $\calO'$-module
to a direct sum $\oplus M'_n$ such that $(\O,M_n)$ is weakly
equivalent to $(\calO',M'_n)$ for each $n$.
\end{prop}
\begin{proof}
It is enough to prove that for  a direct weak equivalence
$$(\Psi,\Phi):(\calO,M)\stackrel{\simeq}{\longrightarrow} (\calR,N),$$
$M$ splits as a direct sum if and only if $N$ splits in a
compatible way.

In one direction, suppose  that $N\simeq\oplus_n N_n$ as
$\calR$-modules. It is clear that the restriction of scalars
functor $\Psi^*$ preserves direct sums and weak equivalences (quasi-isomorphisms).
Therefore $\Psi^*(N)\simeq\oplus_n\Psi^*(N_n)$. Since by
hypothesis $M$ is weakly equivalent to $\Psi^*(N)$, we have the
required splitting of $M$.

In the other direction suppose that the $\calO$-module $M$ is
weakly equivalent to  $\oplus_n M_n$. We can assume that each
$M_n$ is cofibrant, hence so is $\oplus_n M_n$. Moreover
$\Psi^*(N)$ is fibrant because every $\calO$-module is. Therefore,
since $M$ is weakly equivalent to $\Psi^*(N)$, there exists a
direct weak equivalence $\gamma\colon\oplus_n M_n\quism
\Psi^*(N)$. Since $(\Psi^*,\Psi_*)$ is a Quillen equivalence, the
weak equivalence $\gamma$ induces an adjoint weak equivalence
$\gamma^\flat\colon\Psi_*(\oplus_n M_n)\quism N$. As a left adjoint, $\Psi_*$ commutes with
coproducts, therefore we get the splitting
$\oplus_n\Psi_*(M_n)\quism N$. Moreover we have a weak equivalence
$M_n\quism\Psi^*\Psi_*(M_n)$ because it is the adjoint of the
identity map on $\Psi_*(M_n)$, $M_n$ is cofibrant, and
$(\Psi^*,\Psi_*)$ is a Quillen equivalence. Thus that splitting of 
$N$ is compatible with the given splitting of $M$.
\end{proof}


\subsection{Lax monoidal functors, enriched categories, and their
modules}\label{S:lax}


 Let $\calC$ and $\calD$ be two symmetric
monoidal categories. A \emph{lax symmetric monoidal functor}
$F\colon{\mathcal C}\to \mathcal D$ is a (non enriched) functor,
together with  morphisms $\unit_{\calD}\to F(\unit_\calC)$ and
$ F(X)\otimes F(Y)\to F(X\otimes Y)$, natural in
$X,Y\in\calC$, that satisfy the obvious unit, associativity, and
symmetry relations. In this paper, we will some times use ``monoidal'' to mean ``lax symmetric monoidal'', as this is the only notion of monoidality that we will consider.

Such a lax symmetric monoidal functor $F$ induces a functor (which
we will still denote by $F$) from $\CI$-categories to ${\mathcal
D}I$-categories. Explicitly if $\calO$ is a $\calC I$-category
then $F(\calO)$ is the $\calD$-category whose set of objects is
$I$ and morphisms are $(F(\calO))(i,j):=F(\calO(i,j))$. Moreover,
$F$ induces a functor from $\rmods{\O}$ to $\rmods{F(\O)}$. We will denote this
functor by $F$ as well.

The main examples that we will consider
are those from Sections \ref{S:Homology} and \ref{S:Dold-KanCorrespondence}, and their composites:
\begin{enumerate}

\item Homology:
$\HH_*\colon(\Ch,\otimes,\Field)\longrightarrow(\Ch,\otimes,\Field)$;
\item Normalized singular chains:
$\Sing\colon(\Top,\times,*)\longrightarrow(\Ch,\otimes,\Field)\,,\,X\longmapsto
\Sing(X)$.
\end{enumerate}
The fact that the normalized chains functor is lax monoidal, and equivalent to the unnormalized chains functor, is
explained in  \cite[Section 2]{SS}. As is customary, we 
often abbreviate the composite $\HH_*\circ\Sing$ as $\HH_*$.

Recall that we also use the functor $\HH_n\colon
(\Ch,\otimes,\Field)\to(\Ch,\otimes,\Field)$, where $\HH_n(C,d)$ is
seen as a chain complex concentrated in degree $n$. The functor
$\HH_n$ is not monoidal for $n>0$. However, $\HH_0$
is monoidal.

Thus if $\calB$ is a small $\Top\! I$-category then $\Sing(\calB)$ and $\HH_*(\calB)$
are $\ChI$-categories. Also if $B\colon\calB\to\Top$ is a
$\calB$-module then $\Sing(B)$ is a $\Sing(\calB)$-module and
$\HH_*(B)$ is an $\HH_*(\calB)$-module. We  also have the $\ChI$-category
$\HH_0(\calB)$.

\subsection{Discretization of enriched categories}\label{S:discretisation}

When we want to emphasize that a category is
\emph{not} enriched (or, equivalently, enriched over $\Set$), we will use the term \emph{discrete category}. When we speak of an $\calA$-diagram we
always assume that $\calA$ is a discrete category.

Let $\calC$ be a closed symmetric monoidal category. There is a forgetful functor 
$\phi\colon\calC\to \Set$, defined by $$\phi(C):=\hom_{\calC}(\unit, C)$$
It is immediate from the definitions that $\phi$ is a monoidal functor.
Therefore, it induces
a functor from categories enriched over $\calC$ to discrete categories. We will call this induced functor the \emph{discretization} functor. 
Let $\calO$ be a category enriched over $\calC$. The discretization of $\calO$
will be denoted $\calO^\delta$. It has the same objects  as $\calO$, and its sets of morphisms are given by the discretization of morphisms in $\calO$. For example, $\Top$ can be either the $\Top$-enriched category or the associated discrete category.  For $\Ch$, the set of morphisms between two chain complexes $X_*$ and $Y_*$ in the discretization of $\Ch$ is the set of cycles of degree 0 in the chain complex $\Ch(X_*, Y_*)$, i.e. the set of chain maps.  It is easy to see that if $\calC$ is a closed symmetric monoidal category, then the discretization of $\calC$ is the same as $\calC$, considered as a discrete category.  We will not use special  notation to distinguish between $\calC$ and its underlying discrete category.

Let $M\colon\calO\to\calR$ be a $\calC$-functor between two
$\calC$-categories. The \emph{underlying discrete functor} is the
functor
$$M^\delta\colon\calO^\delta\longrightarrow\calR^\delta$$
induced in the obvious way from $M$. More precisely, if $i$ is an
object of $\calO$ then $M^\delta(i)=M(i)$. If $j$ is another
object and $f\in\calO^\delta(i,j)$, that is
$f\colon\unit\to\calO(i,j)$, then $M^\delta(f)\in
\calR^\delta(M^\delta(i),M^\delta(j))$ is defined as the composite
$\unit\stackrel{f}\to\calO(i,j)\stackrel{M{(i,j)}}\longrightarrow\calR(M(i),M(j))$.
Similarly if $\Phi\colon M\to M'$ is an enriched natural
transformation between enriched functors, we have an induced
discrete natural transformation $\Phi^\delta\colon M^\delta\to
M'^\delta$. In particular, an $\calO$-module $M$ induces an
$\calO^\delta$-diagram $M^\delta$ in $\calC$ and a morphism of $\calO$-modules induces a
morphism of $\calO^\delta$-diagrams.

Let $F\colon\calC\to\calD$ be a lax symmetric monoidal functor,
let $\calO$ be a $\calC I$-category, and let $M\colon\calO\to\calC$
be an $\calO$-module. As explained before, we have an induced $\calD I$-category $F(\calO)$, and an  
$F(\calO)$-module $F(M)$.  We may compare $\calO^\delta$ and $F(\calO)^\delta$ 
by means of a functor
$$F^\delta_{\calO}\colon \calO^\delta\longrightarrow F(\calO)^\delta
$$
which is the identity on objects and if
$f\colon\unit_\calC\to\calO(i,j)$ is a morphism in $\calO^\delta$,
then $F^\delta_\calO(f)$ is the composite $\unit_\calD\to
F(\unit_\calC)\stackrel{F(f)}\to F(\calO(i,j))$.

It is straightforward to verify the following two
properties of discretization.
\begin{lemma}\label{L:commuteFdiscrete}
Let $F\colon\calC\to\calD$ be a lax symmetric monoidal functor,
let $\calO$ be a $\calC I$-category and let $M$ be an
$\calO$-module. The following diagram of discrete functors
commutes
$$\xymatrix{
\calO^\delta\ar[d]_-{F^\delta_\calO}\ar[r]^-{M^\delta}&
\calC\ar[d]^-{F}\\
F(\calO)^\delta\ar[r]_-{F(M)^\delta}& \calD. }$$
\end{lemma}

\begin{lemma}
\label{L:discreteweakequiv} Let $\calC$ be a monoidal model
category and  let $\calO$ be a $\calC I$-category. If $\Phi\colon
M\stackrel{\simeq}{\longrightarrow} M'$  is a weak equivalence  of $\calO$-modules then
$\Phi^\delta\colon M^\delta\stackrel{\simeq}{\longrightarrow}  M'^\delta$ is a weak equivalence
of $\calO^\delta$-diagrams.
\end{lemma}
\section{Operads and associated enriched categories}\label{S:operads}

We will first recall the
notions of operads, right modules over operads, and weak
equivalences of operads. We will then describe the enriched category
associated to an operad. Finally, we will treat the central example of
the little balls operad. The 
enriched category viewpoint will help us to deduce (in
\refS{FormalityEmbedding}) the formality of certain
topological functors from the
formality of the little ball operads.

\subsection{Operads and right modules}\label{S:OperadsModules}

Among the many references for operads, a recent one that covers them from a viewpoint similar to ours is Ching's paper \cite{Ching}. However, there is one
important difference between our setting and Ching's:  He only considers operads without the zero-th term, while we consider operads with one. Briefly, 
an \emph{operad} in a symmetric monoidal category
$(\calC,\otimes,\unit)$, or a $\calC$-operad, is a symmetric sequence
$O(\bullet)=\{O(n)\}_{n\in\BbbN}$  of objects of $\calC$, equipped with 
structure maps 
$$O(n)\otimes
O(m_1)\otimes\ldots\otimes O(m_n)\longrightarrow O(m_1+\cdots+m_n)\ \ \text{and }\  
\unit\longrightarrow O(1),
$$ 
satisfying certain associativity, unit, and
symmetry axioms. There is an obvious notion of a \emph{morphism of operads}.

When $\calC$ is a monoidal model category, we say that a morphism
$f\colon O(\bullet)\to R(\bullet)$ of $\calC$-operads is a
\emph{weak equivalence} if $f(n)$ is a weak equivalence in $\calC$
for each natural number $n$. If $f\colon O(\bullet)\to R(\bullet)$
and $f'\colon O'(\bullet)\to R'(\bullet)$ are morphisms of
operads, a \emph{morphism of arrows} from $f$ to $f'$ is a pair
$\left(o\colon O(\bullet)\to O'(\bullet)\,,\,r\colon R(\bullet)\to
R'(\bullet)\right)$ of morphisms of operads such that the obvious
square diagrams commute. Such a pair $(o,r)$ is called a
\emph{weak equivalence} if both $o$ and $r$ are weak
equivalences.

A \emph{right module} over a $\calC$-operad $O(\bullet)$ is a
symmetric sequence $M(\bullet)=\{M(n)\}_{n\in \BbbN}$ of objects of $\calC$, equipped with
structure morphisms
$$M(n)\otimes
O(m_1)\otimes\ldots\otimes O(m_n)\longrightarrow M(m_1+\cdots+m_n)
$$
satisfying certain obvious associativity, unit, and symmetry axioms (see \cite{Ching} for details).
Notice that a morphism of operads $f\colon O(\bullet)\to R(\bullet)$ endows $R(\bullet)$ 
with the structure of a right
$O(\bullet)$-module. 

\subsection{Enriched category associated to an operad}\label{S:OperadEnrichedCategories}

Fix a closed symmetric monoidal category $\calC$ that admits
finite coproducts.
 Recall from Section \ref{S:EnrichedCategories} that a
$\calC\BbbN$-category is a category enriched over $\calC$ whose 
set of objects is $\BbbN$. The
\emph{$\calC\BbbN$-category associated to the $\calC$-operad}
$O(\bullet)$ is the category $\calO$ defined by
$$\calO(m,n)=\coprod_{\alpha\colon\underline{m}\longrightarrow \underline{n}} O({\alpha^{-1}(1)})\otimes \cdots\otimes O({\alpha^{-1}(n)})
$$
where the coproduct is taken over set maps
$\alpha\colon\underline{m}:=\{1,\cdots,m\}\to
\underline{n}:=\{1,\cdots,n\}$ and $O(\alpha^{-1}(j))=O(m_j)$
where $m_j$ is the cardinality of $\alpha^{-1}(j)$. Composition of morphisms
is prescribed by operad structure maps in $O(\bullet)$. In particular $\calO(m,1)=O(m)$.

Let $O(\bullet)$ be a $\calC$-operad  and let $\calO$ be the
associated $\calC\BbbN$-category. A right module (in the sense of operads) $M(\bullet)$ over $O(\bullet)$ gives rise to a right $\calO$-module (in the sense of \refS{enriched}) 
\begin{align*}
M(-)\colon\calO  & \longrightarrow\calC    \\
               n & \longmapsto M(n)
\end{align*}
where $M(-)$ is defined on morphisms by the $\calC$-morphisms
$$M(m,n)\colon\calO(m,n)\longrightarrow\calC(M(n),M(m))$$
obtained by adjunction from the structure maps
$$M(n)\otimes \calO(m,n)=\coprod_{\alpha\colon\underline m\longrightarrow\underline n} M(n)\otimes
O(\alpha^{-1}(1))\otimes\cdots\otimes O(\alpha^{-1}(n))\longrightarrow M(m).$$
 If
$f\colon O(\bullet)\to R(\bullet)$ is a morphism of operads then
 we have an associated right $\calO$-module
$R(-)\colon\calO\to\calC$.

It is obvious that if $O(\bullet)$ and
$O'(\bullet)$ are weakly equivalent, objectwise cofibrant, operads over a monoidal model
category $\calC$ then the associated $\calC\BbbN$-categories
$\calO$ and $\calO'$ are weakly equivalent. Also, if $f\colon
O(\bullet)\to R(\bullet)$ and $f'\colon O'(\bullet)\to
R'(\bullet)$ are weakly equivalent morphisms of operads, then the pair 
$(\calO,R(-))$ is weakly equivalent, in the category of $\calC\BbbN$-modules,  to the
pair $(\calO',R'(-))$.

Let $F\colon\calC\to\calD$ be a lax symmetric monoidal functor, and
suppose $O(\bullet)$ is an operad in $\calC$. Let $\calO$ be
the $\calC\BbbN$-category associated to $O(\bullet)$. Then $F(O(\bullet))$
is an operad in $\calD$, and $F(\calO)$ is a $\calD\BbbN$-category. It is
easy to see that there is a natural morphism from the
$\calD\BbbN$-category associated to the $\calD$-operad
$F(O(\bullet))$ to $F(\calO)$. This morphism is not an isomorphism, unless $F$ is strictly monoidal and also takes coproducts to coproducts, but in all cases that we consider, it will be a weak equivalence. Similarly if $f\colon O(\bullet)\to R(\bullet)$
is a morphism of operads in $\calC$ and if $R(-)$ is the right
$\calO$-module associated to the $O(\bullet)$-module $R(\bullet)$,
then $F(R(-))$ has a natural structure of an $F(\calO)$-module, extending the structure of an 
$F(O(\bullet))$-module possessed by $F(R(\bullet))$.

\subsection{The standard little balls operad}\label{S:LittleDisksOperad}

The most important operad for our purposes is what we will call
the standard balls operad. Let $V$
be a Euclidean space. By a \emph{standard ball} in $V$ we
mean a subset of $V$ that is obtained from the open unit ball by
dilation and translation. The operad of standard balls will be
denoted by $\Balls(\bullet,V)$. It is the well-known operad in $(\Top,\times,*)$, 
consisting of the topological spaces
$$
\Balls(n,V)=\{n\textrm{-tuples of disjoint standard balls inside the unit
ball of } V\}
$$
with the structure maps given by composition of inclusions after
suitable dilations and translations.

The $\Top\!\BbbN$-category associated to the standard balls operad
$\Balls(\bullet,V)$ will be denoted by $\calB(V)$. An object of $\calB(V)$ is 
a non-negative integer 
$n$ which can be thought  of as an abstract (i.e., not embedded) disjoint union
of $n$ copies of the 
unit ball in $V$. The space of morphisms $\calB(V)(m,n)$ is
the space of  embeddings of $m$ unit balls into $n$ unit balls, that on each
ball are obtained by dilations and translations.

Let $j\colon W\hookrightarrow V$ be a linear isometric inclusion of Euclidean spaces. Such a map induces a  morphism of operads
$$
j\colon \Balls(\bullet,W)\longrightarrow \Balls(\bullet,V)
$$
where a ball centered at $w\in W$ is sent to the ball of the same radius centered
at $j(w)$.

Hence $\Balls(\bullet,V)$ is a right module over
$\Balls(\bullet,W)$, and we get a right $\calB(W)$-module
\begin{align*}
\Balls(-,V)\colon \calB(W) & \longrightarrow\Top     \\
                         n & \longmapsto \Balls(n,V).
\end{align*}

We can  apply lax monoidal functors to the above setting. For
example, $\Sing(\Balls(\bullet,W))$ and $\HH_*(\Balls(\bullet,W))$
are operads in $(\Ch,\otimes,\Field)$. Hence we get
$\Ch\!\BbbN$-categories $\Sing(\calB(W))$ and $\HH_*(\calB(W))$, a
right  $\Sing(\calB(W))$-module $\Sing(\Balls(-,V))$, and a right
$\HH_*(\calB(W))$-module $\HH_*(\Balls(-,V))$.

We will also consider the discrete categories $\calB(W)^\delta$ and $\Sing(\calB(W))^\delta$
obtained by by the discretization process from $\calB(W)$ and $\Sing(\calB(W))$ respectively. Note that $\Sing(\calB(W))^\delta=\Field[\calB(W)^\delta]$.

\section{Formality and splitting of the little balls
operad}\label{S:formalityballs}
In this section, all chain complexes and homology groups are taken with coefficients in $\R$. 
A deep theorem of Kontsevich  (\refT{KontsevichFormality} of the
Introduction and Theorem 2 of \cite{K_Formality}) asserts that the standard balls operad is formal over the reals. We will
need a slight strengthening of this result.  Throughout this section, let $j\colon W\hookrightarrow V$
 be, as usual, a linear isometric inclusion of Euclidean spaces.
Recall the little balls operad and the
associated enriched categories and modules as in
\refS{LittleDisksOperad}. Here is the version of Kontsevich's theorem we need.
\begin{thm}[{\bf Relative Formality}]\label{T:relativeKontsevich}
 If $\dim V>2\dim W$ then the morphism of chain operads
 $$\Sing(j)\colon \Sing(\Balls(\bullet,W))\otimes \R\longrightarrow
\Sing(\Balls(\bullet,V))\otimes \R$$ is weakly equivalent to  the morphism
$$\HH_*(j)\colon \HH_*(\Balls(\bullet,W); \R)\longrightarrow
\HH_*(\Balls(\bullet,V);\R).$$
\end{thm}
\begin{proof}[Sketch of the proof]
\newcommand{\FM}{\operatorname{FM}}
\newcommand{\SemiAlgChains}{\textsf{SemiAlgChain}}
\newcommand{\Graphs}{\textsf{Graphs}}
A detailed proof will appear in \cite{LV_Formality}. Here we give a sketch based on the proof absolute formality given in \cite[Theorem 2]{K_Formality}, and we follow that paper's notation.
Denote by $\FM_d(n)$ the Fulton-MacPherson compactification of the configuration space of $n$ points
in $\BbbR^d$. This defines an operad $\FM_d(\bullet)$ which is homotopy equivalent to the little balls operad $B(\bullet,\BbbR^d)$.
Kontsevich constructs a quasi-isomorphism
$$
\Psi\colon\SemiAlgChains_*(\FM_d(n))\stackrel{\simeq}{\longrightarrow} \Graphs_d(n)\hat\otimes\BbbR
$$
where $\SemiAlgChains_*$ is a chain complex of semi-algebraic chains naturally quasi-isomorphic to
singular chains and $\Graphs_d$ is the chain complex of admissible graphs defined in \cite[Definition 13]{K_Formality}.
For $\xi$ a semi-algebraic chain on $\FM_d(n)$, the map $\Psi$ is defined by
$$\Psi(\xi)=\sum\Gamma\otimes\langle\omega_\Gamma,\xi\rangle,$$
where the sum is taken over all admissible graphs $\Gamma$ and $\omega_\Gamma$ is the differential form defined
in \cite[Definition 14]{K_Formality}.

Let $j_*\colon \FM_{\dim W}(n)\to\FM_{\dim V}(n)$ be the map induced by the inclusion of Euclidean spaces $j$. Notice that $H_i(j_*)=0$
for $i>0$. Define $\epsilon\colon  \Graphs_{\dim W}(n)\to\Graphs_{\dim V}(n)$ to be zero on graphs with at least one edge, and the identity
on the graph without edges.
We need to show that the following diagram commutes:

$$\xymatrix{\SemiAlgChains_*(\FM_{\dim W}(n))\ar[d]^-{j_*}
\ar[r]^-\simeq& \Graphs_{\dim W}(n)\hat\otimes\BbbR\ar[r]^-\simeq\ar[d]^-\epsilon&\ar[d]^-{\HH(j_*)}\HH_*(\FM_{\dim W}(n))\\
\SemiAlgChains_*(\FM_{\dim V}(n))\ar[r]^-\simeq &\Graphs_{\dim V}(n)\hat\otimes\BbbR\ar[r]^-\simeq&\HH_*(\FM_{\dim V}(n)).
}
$$
The commutativity of the right hand square is clear. For the left hand square it suffices to check that for any admissible graph of positive degree
$\Gamma$
and for any non-zero semi-algebraic chain $\xi\in \SemiAlgChains_*(\FM_{\dim W}(n))$ we have $\langle\omega_\Gamma,j_*(\xi)\rangle=0$.

The first $n$ vertices of $\Gamma$,  $1,\cdots,n$, are called external and the other are called internal.
If every external vertex of $\Gamma$ is connected to an edge, then, using the fact that internal vertices are at least trivalent, we obtain that
the form $\omega_\Gamma$ on $\FM_{\dim V}(n)$ is of degree $\geq n(\dim V-1)/2$. Since $\dim V >2\dim W$, we get that $\deg(\omega_\Gamma)>\dim (\FM_{\dim W}(n))$.  Therefore $\deg(\omega_\Gamma)>\deg(j_*(\xi))$ and $\langle\omega_\Gamma,j_*(\xi)\rangle=0$.

If $\Gamma$ has an isolated external vertex, then $\langle\omega_\Gamma,j_*(\xi)\rangle=\langle\omega_\Gamma,j_*(\xi')\rangle$,
where $\xi'$ is a chain in $\SemiAlgChains_*(\FM_{\dim W}(m))$ with $m<n$ and the proof proceeds by induction.
\end{proof}
We remark once again that the formality theorem is for chain complexes over $\R$, not over $\Q$. We do not know if the little balls operad is formal over the rational numbers, but we do think it is an interesting question. We note that a general result about descent of formality from $\R$ to $\Q$ was proved in \cite{GNPR}, for operads \emph{without a term in degree zero}. The proof does not seem to be easily adaptable to operads with a zero term.

To deduce the formality of certain diagrams more directly related to spaces of embeddings, 
we first reformulate relative formality in terms of homogeneous
splittings in the spirit of \refP{PostnikovInterpretation}. With
this in mind we introduce the following enrichment of
\refD{discrete-splits}.
\begin{definition}\label{D:module-splits}
Let $\calO$ be a $\ChI$-category. We say that an $\calO$-module
$M\colon\calO\to\Ch$  \emph{splits homogeneously} if there exists
a sequence $\{ M_n\}_{n\in\BbbN}$ of $\calO$-modules such that
$M\simeq \oplus_n M_n$ and $\HH_*(M_n)=\HH_n(M_n)$.
\end{definition}
Our first example (a trivial one) of such a  homogeneous splitting of modules is
given by the following
\begin{lemma}\label{L:H*Bsplits}
If $\dim V>\dim W$ then the $\HH_*(\calB(W))$-module
$\HH_*(\Balls(-,V))$ splits homogenously.
\end{lemma}
\begin{proof}
Notice that $\HH_0(\calB(W))$ is also a $\Ch\!\BbbN$-category.
Since  our chain complexes are non-negatively graded and with a zero differential, 
 we
have an obvious inclusion functor
$$i\colon  \HH_0(\calB(W))\hookrightarrow  \HH_*(\calB(W))$$
 and a projection functor
$$\Phi\colon \HH_*(\calB(W))\longrightarrow \HH_0(\calB(W))$$
between $\Ch\!\BbbN$-categories,
where $\Phi\circ i$ is the identity.
Therefore, an $\HH_*(\calB(W))$-module admits a structure of an
$\HH_0(\calB(W))$-module via $i$. Since
$\HH_0(\calB(W))$ is a category of chain complexes concentrated in degree $0$ and
$\HH_*(\Balls(-,V))$ has no differentials, it is clear that we have
a splitting of $\HH_0(\calB(W))$-modules
\begin{equation}\label{E:splitH0}
\HH_*(\Balls(-,V))\cong\oplus_{n=0}^\infty \HH_n(\Balls(-,V)).
\end{equation}

Moreover, since $\dim W<\dim V$ the morphisms
$$
\HH_*(\Balls(n,W))\longrightarrow \HH_*(\Balls(n,V))
$$
are zero in positive degrees. Hence the $\HH_*(\calB(W))$-module
structure on $\HH_*(\Balls(-,V))$ factors through
the above-mentioned $\HH_0(\calB(W))$-module structure via $\Phi$. Therefore, the splitting  $(\ref{E:splitH0})$ is a splitting of $\HH_*(\calB(W))$-modules.
\end{proof}
Using \refL{H*Bsplits} and the Relative Formality Theorem, we obtain the following highly non-trivial splitting.
\begin{lemma}\label{L:C*Bsplits}
If  $\dim V>2\dim W$ then the $\Sing(\calB(W))$-module
$\Sing(\Balls(-,V))$ splits homogenously.
\end{lemma}
\begin{proof}
We deduce from  \refT{relativeKontsevich} that the $\Ch\!\mathbb{N}$-module categories
$(\Sing(\calB(W)), \Sing(\Balls(-,V)))$ and
 $(\HH_*(\calB(W)), \HH_*(\Balls(-,V)))$ are equivalent. By \refL{H*Bsplits} the latter splits
homogeneously, hence, by \refP{SplittingPropagates} the same is
true of the former.
\end{proof}

Recall from \refS{discretisation} that the enriched category
$\calB(W)$ has an underlying discrete category $\calB(W)^\delta$
and that the $\calB(W)$-module $\Balls(-,V)$ induces a
$\calB(W)^\delta$-diagram $\Balls(-,V)^\delta$.
\begin{prop}\label{P:formalitydiscretedisks}
If $\dim V>2\dim W$ then the $\calB(W)^\delta$-diagram
$$\Sing(\Balls(-,V))^\delta\colon\calB(W)^\delta\longrightarrow\CC_{\R}$$
is formal.
\end{prop}
\begin{proof}
By  \refL{commuteFdiscrete} the following diagram of discrete
functors commutes:
$$
\xymatrix{
\calB(W)^\delta\ar[d]_-{(\Sing)^\delta_{\calB(W)}}\ar[rr]^-{\Balls(-,V)^\delta} &  &  \Top\ar[d]^{\Sing}\\
\left(\Sing(\calB(W))\right)^\delta\ar[rr]^-{(\Sing(\Balls(V,-)))^\delta}& &    \Ch_\R
}
$$
where the discretization  $\left(\Sing(\calB(W))\right)^\delta$ is nothing else than the linearization 
of $(\calB(W))^\delta$, that is the category whose morphisms are formal $\R$-linear combinations of morphism 
in  $\calB(W)$ and the objects are the same as in $\calB(W)$. 
We want to prove that the $\calB(W)^\delta$-diagram
$\Sing\left(\Balls(-,V)^\delta\right)$ is formal. By the
commutativity of the square above and \refL{formalprecomposite} it is
enough to prove that the $(\Sing(\calB(W)))^\delta$-diagram
$(\Sing(\Balls(V,-)))^\delta$ is formal. By \refL{C*Bsplits} the
$\Sing(\calB(W))$-module $\Sing(\Balls(-,V))$ splits
homogeneously. By Lemma 
\ref{L:discreteweakequiv} we deduce that the
$\Sing(\calB(W))^\delta$-diagram $\Sing(\Balls(-,V))^\delta$
splits homogeneously, which implies by
\refP{PostnikovInterpretation} the formality of that diagram.
\end{proof}


\section{Formality of a certain diagram arising from embedding calculus}
\label{S:FormalityEmbedding}

In this section, all chain complexes are still taken over the real numbers. 
As before, fix a linear isometric inclusion of Euclidean vector
spaces $j\colon W\hookrightarrow V$. Let $\calO(W)$ be the poset of open subsets of $W$. As
explained in the Introduction, we have two contravariant functors 
$$
\Emb(-,V)\,,\,\Imm(-,V)\colon\calO(W)\longrightarrow\Top.
$$
Moreover, the fixed
embedding $j\colon W\hookrightarrow V$ can serve as a basepoint,
so we can consider the homotopy fiber of the inclusion $\Emb(-,V)\to\Imm(-,V)$, which we denote by
$$\Ebar(-,V)\colon\calO(W)\longrightarrow\Top.$$

Our goal in this section is to compare a certain variation of this
functor with the functor
$$
\Balls(-,V)^\delta\colon \calB(W)^\delta\longrightarrow\Top
$$
and to deduce in \refT{formalityembeddingdiagram} the stable
formality of certain diagrams of embedding spaces. In order to do
this we first introduce a subcategory $\Os(W)$ of $\calO(W)$ and a
category $\tOs(W)$ which will serve as a turning table between
$\Os(W)$ and $\calB(W)^\delta$.

To describe $\Os(W)$ recall that a \emph{standard ball} in $W$ is an
open ball in the metric space $W$, i.e. it is
obtained in a unique way by a dilation and translation of the unit ball in $W$.  The category
$\Os(W)$ is the full subcategory of $\calO(W)$ whose objects are
finite unions of disjoint standard balls.

The category $\tOs(W)$ is a kind of covering
of $\Os(W)$. Recall that the object $m\in\BbbN$ of $\calB(W)$ can be thought of as an abstract disjoint union of $m$ copies of the unit ball of
$W$. An object of $\tOs(W)$ is then an embedding $\phi\colon
m\hookrightarrow W$ such that the restriction of $\phi$ to each unit ball amounts to a 
dilation and translation. In other words an object $(\phi,m)$ of
$\tOs(W)$ is the same as an ordered $m$-tuple of disjoint standard balls
in $W$. The union of these $m$ standard balls is an object of
$\Os(W)$ that we denote by $\phi(m)$, as the image of the
embedding $\phi$. By definition, there is a morphism in $\tOs(W)$
between two objects $(\phi,m)$ and $(\psi,n)$   if and only if
$\phi(n)\subset\psi(m)$, and such a morphism is unique.

We define functors
$$
\xymatrix{\calB(W)^\delta&\ar[l]_\lambda
\tOs(W)\ar[r]^\pi&\Os(W)}.
$$
Here $\pi$ is defined on objects by $\pi(\phi,m)=\phi(m)$
and is defined on morphisms by sending a morphism $\alpha:(\phi_1,m_1)\to (\phi_2,m)$ to the inclusion $\phi_1(m_1)\hookrightarrow \phi_2(m_2)$, and this functor is easily seen to be an equivalence of categories.
 The functor $\lambda$
is defined on objects by $\lambda(\phi,m)=m$, and is defined on morphisms using the fact that any two standard balls in $W$ can be canonically identified by a unique transformation that is a combination of dilation and translation. 

We would like to compare the following two composed  functors
\begin{eqnarray*}
\xymatrix{
\Ebar(\pi(-),V))\colon\tOs(W) \ar[r]^-{\pi}  &  \Os(W) \ar[rr]^-{\Ebar(-,V)} &  & \Top } \\
\xymatrix{\Balls(\lambda(-),V)^\delta\colon\tOs(W) \ar[r]^-{\lambda} & \calB(W)^\delta\ar[rr]^-{\Balls(-,V)^\delta}  &  & \Top.
}
\end{eqnarray*}

\begin{prop}\label{P:compareEmbpiBLambda}
The $\tOs(W)$-diagrams $\Balls(\lambda(-),V)^\delta$ and
$\Ebar(\pi(-),V)$ are weakly equivalent.
\end{prop}
\begin{proof}
Define subspaces $\LEmb(\phi(n),V)\subset \Emb(\phi(n),V)$ and
$\LImm(\phi(n),V)\subset \Imm(\phi(n),V)$ to be the spaces of
embeddings and immersions, respectively, that are  affine on each ball. It is well-known that
the above inclusion maps are homotopy equivalences.
We may define $\LEbar(\phi(n),V)$ to be the homotopy
fiber of the map $\LEmb(\phi(n),V) \to \LImm(\phi(n),V)$. Thus
there is a natural homotopy equivalence
$$\LEbar(\phi(n),V)\stackrel{\simeq}{\longrightarrow} \Ebar(\phi(n),V).$$ Define
$\widetilde{\operatorname{Inj}}(W,V)$ as the space of injective linear
maps from $W$ to $V$, quotiented out by the multiplicative group
of positive reals, i.e. defined up to dilation. Then there is a
natural homotopy equivalence
$$\LImm(\phi(n),V)\stackrel{\simeq}{\longrightarrow} \widetilde{\operatorname{Inj}}(W,V)^n$$
obtained  by differentiating the immersion at each component of
$\phi(n)$. Now consider the map 
$$\LEmb(\phi(n),V)\longrightarrow\widetilde{\operatorname{Inj}}(W,V)^n.$$ We denote the homotopy fiber
of this map by $F(n,\phi)$, and we obtain a natural equivalence
$$\LEbar(\phi(n),V)\stackrel{\simeq}{\longrightarrow}  F(n,\phi).$$ Finally since the composite map
$$\Balls(n,V)\hookrightarrow \LEmb(\phi(n),V) \longrightarrow\widetilde{\operatorname{Inj}}(W,V)^n$$ is the constant map into the basepoint,
there is a natural map $\Balls(n,V)\to F(n,\phi)$. It is easy to
see that the map is an equivalence. To summarize, we have
constructed the following chain of natural weak equivalences
$$\Ebar(\phi(n),V) \stackrel{\simeq}{\leftarrow} \LEbar(\phi(n),V) \stackrel{\simeq}{\longrightarrow} F(n,\phi) \stackrel{\simeq}{\leftarrow} \Balls(n,V).$$
\end{proof}
We are ready to prove the main result of this section.
\begin{thm}\label{T:formalityembeddingdiagram}
If $\dim V>2\dim W$ then the $\tOs(W)$-diagram 
$\Sing(\Ebar(\pi(-),V))$ is stably formal.
\end{thm}
\begin{proof}
By \refP{formalitydiscretedisks} and \refL{formalprecomposite} the
diagram $\Sing(\Balls(\lambda(-),V))^\delta$ is stably formal.
Proposition \ref{P:compareEmbpiBLambda} implies the theorem.
\end{proof}

\section{More generalities on calculus of functors}\label{S:ReviewCalculi}

In this section we digress to review in a little more detail the basics of embedding and orthogonal calculus. We will also record some general observations about bi-functors to which both brands of calculus apply. The
standard references are
\cite{WeissEmb} and \cite{WeissOrth}.

\subsection{Embedding calculus} Let $M$ be a  smooth manifold (for convenience, we assume that $M$ is the interior of a compact manifold with boundary). Let $\calO(M)$ be the poset of open subsets of $M$ and let 
$\calO_k(M)$ be the subposet consisting of open subsets homeomorphic to a union of
at most $k$ open balls. Embedding calculus is concerned with the study of contravariant isotopy functors from  $F$ to a Quillen model category 
(Weiss only considers functors into the category of spaces, and, implicitly, spectra, but much of the theory works just as well in the more general setting of model categories). 
Following \cite[page 5]{WeissEmb}, we say that a contravariant functor is \emph{good} if it converts isotopy equivalences to weak equivalences and filtered unions to homotopy limits. 
Polynomial functors are defined in terms of certain cubical diagrams, similarly to the way they are defined in Goodwillie's homotopy calculus. 
Recall that a cubical diagram of spaces is called {\it strongly co-cartesian} if each of its two-dimensional faces is a homotopy pushout square. A contravariant functor 
$F$ on $\calO(M)$ is called {\it polynomial of degree $k$} if it takes strongly co-cartesian $k+1$-dimensional cubical diagrams of opens subsets of $M$ to
homotopy cartesian cubical diagrams (homotopy cartesian cubical diagrams is synonymous with homotopy pullback cubical diagrams).
Good functors can be
approximated by the stages of the tower
defined by
$$
T_kF(U)=\underset{\{U'\in\calO_k(M)\mid U'\subset U\}}\holim\,F(U').
$$

It turns out that $T_kF$ is polynomial of degree $k$, and moreover there is a natural
map $F\longrightarrow T_kF$ which in some sense is the best possible approximation of $F$ by a polynomial functor of degree $k$. More precisely, the map $F\longrightarrow T_kF$ can be characterized as the essentially unique map from $F$ to a polynomial functor of degree $k$ that induces a weak equivalence when evaluated on an object of $\calO_k(M)$. In the terminology of \cite{WeissEmb}, $T_kF$ is the $k$-th Taylor polynomial of $F$. $F$ is said to be {\it homogeneous of degree $k$} if it is polynomial of degree $k$ and $T_{k-1}F$ is equivalent to the trivial functor. 
For each $k$, there is a natural map $T_{k}F\to T_{k-1}F$, compatible with the  maps $F\to T_kF$ and $F\to T_{k-1}F$. Its homotopy fiber is a homogeneous functor of degree $k$, and it  
is called the \emph{$k$-th layer} of the tower.  It plays the role of the $k$-th term in the Taylor series of a function. For functors with values in pointed spaces, 
there is a useful general formula for the $k$-th layer in terms of spaces of
sections of a certain bundle $p:E\to {M \choose k}$ over the space
$M\choose k$ of unordered $k$-tuples of distinct points in $M$. The
fiber of $p$ at a  point $\underline{m}=\{m_1,\ldots,m_k\}$ is
$\widehat{F(m)} $, which is defined to be the total
fiber of the $k$-dimensional cube $S\mapsto F(N(S))$ where $S$
ranges over subsets of $\underline m$ and $N(S)$ stands for a
``small tubular neighborhood'' of $S$ in $M$, i.e., a disjoint
union of open balls in $M$. The fibration $p$ has a preferred
section. See~\cite{WeissEmb}, especially Sections 8 and 9, for
more details and a proof of the following proposition.
\begin{prop}[Weiss]\label{P:LayersOfEmbeddingTower}
The homotopy fiber of the map $T_kF\to T_{k-1}F$ is equivalent to
the space of sections of the fibration $p$ above which agree with
the preferred section in a neighborhood of the fat diagonal in
$M^k$.
\end{prop}
We denote this space of restricted sections by
$$\Gamma_c\left({M \choose k}, \widehat{F(k)}\right).$$

Even though $T_kF$ is defined as the homotopy limit of an infinite category, for most 
moral and practical purposes it behaves as if it was the homotopy limit of a \emph{very small category} (i.e., a category whose simplicial nerve has finitely many non-degenerate simplices). This is so because of the following proposition.
\begin{prop}\label{P:EssentiallyFinite}
There is a very small subcategory $\calC$ of $\calO_k(M)$ such that restriction from $\calO_k(M)$ to $\calC$ induces an equivalence on homotopy limits of all good 
functors.  
\end{prop}
\begin{proof}
It is not difficult to show, using handlebody decomposition and induction (the argument is essentially contained in the proof of Theorem 5.1 of \cite{WeissEmb}) that one can find a finite collection $\{U_1,\ldots,U_N\}$ of open subsets of $M$ such that all their possible intersections are objects of $\calO_k(M)$ and 
$$M^k=\cup_{i=1}^N U_i^k$$
This is equivalent to saying that the sets $U_i$ cover $M$ in what Weiss calls the Grothendieck topology ${\mathcal J}_k$. By \cite{WeissEmb}, Theorem 5.2, polynomial 
functors of degree $k$ are homotopy sheaves with respect to ${\mathcal J}_k$. In practice, this means the following. Let $\calC$ be the subposet of $\calO_k(M)$ given by the sets $U_i$ and all their possible intersections (clearly, $\calC$ is a very small category). Let $G$ be a polynomial functor of degree $k$. Then the following canonical map is a homotopy equivalence
$$G(M)\longrightarrow \underset{U\in \calC}{\holim}\, G(U).$$

We conclude that for a good  functor $F$, there is the following zig-zag of weak equivalences.
$$\underset{U\in\calC}{\holim }\, F(U) \stackrel{\simeq}{\longrightarrow} \underset{U\in\calC}{\holim }\, T_kF (U)\stackrel{\simeq}{\longleftarrow} T_kF(M)$$
Here the left map is a weak equivalence because the map $F\to T_kF$ is a weak equivalence on objects of $\calO_k(M)$, and all objects of $\calC$ are objects of $\calO_k(M)$. The right map is an equivalence because $T_kF$ is a polynomial functor of degree $k$, in view of the discussion above.
\end{proof}
The important consequence of the proposition is that $T_kF$ commutes, up to a zig-zag of weak equivalences, with filtered homotopy colimits of functors. In the same spirit, we have the following proposition.
\begin{prop}\label{P:T_kCommutesWithReals}
Let $F:\calO_k(M)\to \CC_{\Q}$ be a good  functor into rational chain complexes. Then the natural map
$$(T_kF(M))\otimes \R\longrightarrow T_k(F\otimes \R)(M)$$ is a weak equivalence.
\end{prop}
\begin{proof}
Tensoring with $\R$ obviously commutes up to homotopy with very small homotopy limits, and so the claim follows from \refP{EssentiallyFinite}.
\end{proof}

\subsection{Orthogonal calculus} The basic reference for Orthogonal calculus is \cite{WeissOrth}. Let $\calJ$ be the topological category of Euclidean spaces and
linear isometric inclusions. Orthogonal calculus is concerned with the study of \emph{continuous}
functors from $\calJ$ to a model category enriched over $\Top_*$. We will only consider functors into $\Top_*$, $\Spectra$ and closely related categories. 
Like embedding calculus, orthogonal calculus comes equipped with a notion of a polynomial functor, and with a construction that associates with a functor $G$ a tower of approximating functors $P_nG$ such that $P_nG$ is, in a suitable sense, the best possible approximation of $G$ by a polynomial functor of degree $n$. $P_n$ is defined as a certain filtered homotopy colimit of compact homotopy limits. For each $n$, there is a natural map
$P_{n}G\to P_{n-1}G$ and its fiber (again called the \emph{$n$-th layer}) is denoted by $D_nG$. $D_nG$ is a 
homogeneous functor, in the sense that it is polynomial of degree $n$ and $P_{n-1}D_nG\simeq *$.  The following
characterization of homogeneous functors is proved in
\cite{WeissOrth}.
\begin{thm}[Weiss]
Every homogeneous functor of degree $n$ from vector spaces to
spectra is equivalent to a functor of the form
$$\left(C_n\wedge S^{nV}\right)_{h\OO(n)}$$
where $C_n$ is a spectrum with an action of the orthogonal group
$\OO(n)$, $S^{nV}$ is the one-point compactification of the vector
space $\R^n\otimes V$, and the subscript $h\OO(n)$ denotes
homotopy orbits.
\end{thm}
It follows, in particular, that given a (spectrum-valued) functor
$G$ to which orthogonal calculus applies, $D_nG$ has the form
described in the theorem, with some spectrum $C_n$. The spectrum
$C_n$ is called the \emph{$n$-th derivative} of $G$. There is a useful description of
the derivatives of $G$ as stabilizations of certain types of iterated cross-effects of $G$. 

Let $G_1,G_2$ be two functors to which orthogonal calculus applies. Let $\alpha:G_1\to G_2$ be a natural transformation. Very much in the spirit of Goodwillie's homotopy calculus, we say that $G_1$ and $G_2$ {\it agree to $n$-th order via $\alpha$} if the map $\alpha(V):G_1(V)\to G_2(V)$ is $(n+1)\dim(V) +c$-connected, where $c$ is a possibly negative constant, independent of $V$. Using the description of derivatives in terms of cross-effects, it is easy to prove the following proposition
\begin{prop}\label{P:ConnectivityEstimates}
Suppose that $G_1$ and $G_2$ agree to $n$-th order via a natural transformation $\alpha\colon G_1\to G_2$. Then $\alpha$ induces an equivalence on the first $n$ derivatives, and therefore an equivalence on $n$-th Taylor polynomials
$$P_n\alpha:P_nG_1\stackrel{\simeq}{\longrightarrow} P_nG_2$$
\end{prop}

\subsection{Bifunctors} In this paper we consider bifunctors
$$
E\colon\calO(M)^{\operatorname{op}}\times\calJ\longrightarrow\Top/\Spectra
$$
such that the adjoint contravariant functor $\calO(M)\to {\operatorname{Funct}}(\calJ,\Top/\Spectra)$ is good (in the evident sense) and the adjoint functor $\calJ\to\operatorname{Funct}(\calO(M)^{\operatorname{op}},\Top/\Spectra)$ is continuous.  We may apply both embedding calculus and orthogonal calculus to such a bifunctor. Thus by $P_nE(M,V)$ we mean the functor obtained from $E$ by considering it a functor of $V$, (with $M$ being a ``parameter'') and taking the $n$-th Taylor polynomial in the orthogonal sense. Similarly, $T_kE(M,V)$ is the functor obtained by taking the $k$-th Taylor polynomial in the sense of embedding calculus.

We will need a result about the interchangeability of order of applying the differential operators $P_n$ and $T_k$. Operator $T_k$ is constructed using a homotopy limit, while $P_n$ is constructed using a homotopy limit (over a compact topological category) and a filtered homotopy colimit. It follows that there is a natural transformation 
$$P_nT_kE(M,V) \longrightarrow T_kP_nE(M,V)$$
and a similar natural transformation where $P_n$ is replaced with $D_n$.
\begin{lemma}\label{L:PnAndTkCommute} Let $E$ be a bifunctor as above. 
For all $n$ and $k$ the natural map 
$$P_n T_k E(M,V) \stackrel{\simeq}{\longrightarrow} T_kP_nE(M,V)$$
is an equivalence. There is a similar equivalence where $P_n$ is replaced by $D_n$.
\end{lemma}
\begin{proof}
By \refP{EssentiallyFinite}, $T_k$ can be presented as a very small homotopy limit. Therefore, it commutes up to homotopy with homotopy limits and filtered homotopy colimits. $P_n$ is constructed using homotopy limits and filtered homotopy colimits. Therefore, $T_k$ and $P_n$ commute.
\end{proof}


\section{Formality and the embedding
tower}\label{S:RatEmb}


In this section we assume that $\alpha\colon
M\hookrightarrow W$ is an inclusion of an \emph{open subset} into a Euclidean space $W$. From our point of view, there is no loss of generality in this assumption, because if $M$ is an embedded manifold in $W$, we can replace $M$ with an open tubular neighborhood, without changing the homotopy type of $\Ebar(M,V)$. 
As usual, we fix an
isometric inclusion $j\colon W\hookrightarrow V$ of Euclidean vector
spaces.

Recall that we defined the functor
$$\Ebar(-,V)\colon\calO(M)\longrightarrow\Top.$$
The stable rationalisation $\HQ\wedge\Ebar(-,V)_+$ of $\Ebar(-,V)$
admits a Taylor tower (in this section, Taylor towers are taken in the sense of
embedding calculus). Our goal is to give in
 \refT{LimPostnikovSplits} a splitting of the $k$-th stage of this tower. The splitting is
 \emph{not} as a product of the layers in the embedding towers. Rather, we will see in the next section that the splitting is as a product of the layers in the orthogonal tower.

Recall the poset $\Os(W)$ of finite unions of standard balls in
$W$ from \refS{FormalityEmbedding}. Let $\Os(M)$ be the full
subcategory of $\Os(W)$ consisting of the objects which are
subsets of $M$.  For a natural number 
$k$ we define $\Os_k(M)$ as the full subcategory of $\Os(M)$
consisting of disjoint unions of at most $k$ standard balls in $M$.

\begin{prop}\label{P:Tk-tOs}
Let $M$ be an open submanifold of a vector space $W$ and let
$F\colon\calO(M)\to\Top$ be a good functor. The restriction map
$$T_kF(M):= \underset{U\in\calO_k(M)}\holim F(U)\longrightarrow \underset{U\in \StrictOO_k(M)}{\holim} F(U),$$
induced by the inclusion of categories $\StrictOO_k(M)\to \calO_k(M)$, is a homotopy equivalence.
\end{prop}
\begin{proof}
Define 
$T^s_kF(M):=\underset{U\in\StrictOO_k(M)}\holim F(U)$.  There are projection maps 
$$T^s_kF(M)\longrightarrow T^s_{k-1}F(M)$$
induced by the inclusion of categories  $\StrictOO_{k-1}(M)\to\StrictOO_k(M)$, and the map $T_kF\to T^s_kF$ extends to a map of towers. One can adapt the methods of \cite{WeissEmb} to analyze the functors $T_k^sF$. In particular, it is not hard to show, using the same methods as in \cite{WeissEmb}, that our map induces a homotopy equivalence from the homotopy fibers of the map $T_kF\to T_{k-1}F$ to the homotopy fibers of the map $T_k^sF\to T_{k-1}^sF$, for all $k$. Our assertion follows by induction on $k$. 
\end{proof}

Recall the category $\tOs(W)$ defined in \refS{FormalityEmbedding}. Let $\tOs(M)$ be the full subcategory of $\tOs(W)$ consisting of
objects $(\phi,m)$ such that $\phi(m)$ is a subset of $M$. Define
also $\tOsk(M)$ to be the full subcategory of $\tOs(W)$ consisting
of objects $(\phi,m)$ such that $m$ is at most $k$.

Recall the functor
$\pi\colon\tOs(W)\to\Os(W),(\phi,m)\mapsto\phi(m)$, defined in \refS{FormalityEmbedding}. It is clear that this functor restricts to a
functor $\pi\colon\tOsk(M)\to\Os_k(M)$. Moreover it is an equivalence of categories, therefore pullbacks along $\pi$ induce weak equivalences
between homotopy limits.  

We can now prove the main result of this section. Recall
from \refS{FormalityEmbedding} the functor
$$\Balls(\lambda(-),V)\colon\tOs(W)\to\Top$$ which by abuse
of notation we denote by $(\phi,m)\mapsto B(m,V)$.


\begin{thm}\label{T:LimPostnikovSplits}
Let $W\subset V$ be an inclusion of Euclidean vector spaces, let
$M$ be an open submanifold of  $W$, and let $k$ be a natural
number. If $\dim V>2\dim W$ then there is an equivalence of
spectra
$$T_k \HQ\wedge \Ebar(M,V)_+ \simeq \prod _{i=0}^\infty T_k ||\HH_i \Ebar(M,V)||\simeq\prod_{i=0}^\infty \underset{(\phi,m)\in \tOsk(M)}{\holim} ||\HH_i(\Ebar(\pi(\phi,m),V)))|| $$
where $||\HH_i(X)|| $ is the Eilenberg-Mac Lane spectrum
that has the $i$-th rational homology of $X$ in degree
$i$.
\end{thm}
\begin{proof}
By \refP{Tk-tOs} and since $\pi$ is an equivalence of categories we have
$$
 T_k \HQ\wedge \Ebar(M,V)_+ \simeq\underset{(\phi,m)\in\tOsk(M)}\holim\HQ\wedge \Ebar(\pi(\phi,m),V))_+.
$$
By \refP{compareEmbpiBLambda}, the functors 
$\Ebar(\pi(\phi,m),V)$ and $\Balls(\lambda(\phi,m),V)=\Balls(m,V)$ are weakly equivalent, as functors
on $\tOsk(W)$. It follows that their restrictions to $\tOsk(M)$ are weakly equivalent, and so
$$T_k\HQ\wedge \Ebar(M,V)_+\simeq\underset{(\phi,m)\in\tOsk(M)}\holim\HQ\wedge \Balls(m,V)_+.$$
Using the Quillen equivalence between rational spectra and rational chain complexes, and the fact that homotopy limits are preserved by Quillen equivalences, we conclude that 
there is a weak equivalence (or more precisely a zig-zag of weak equivalences) in $\CC_{\Q}$
$$T_k\Sing (\Ebar(M,V))\simeq\underset{(\phi,m)\in\tOsk(M)}\holim\, \Sing(\Balls(m,V)).$$
On the other hand, by \refP{formalitydiscretedisks} and \refL{formalprecomposite}, the functor $m\mapsto \Sing(\Balls(m,V))\otimes \R$ from $\tOsk(M)$ to $\CC_{\R}$ is formal. By 
Remark \ref{R:FormalLimitSplits} we get that
$$
\underset{(\phi,m)\in\tOsk(M)}\holim\Sing(\Ebar(\pi(\phi,m),V))\otimes \R \simeq \prod_{i=0}^\infty
\underset{(\phi,m)\in \tOsk(M)}{\holim}
\HH_i(\Balls(m,V);\R).$$
Recall that $\Balls(m,V)$ is equivalent to the space of configurations of $m$ points in $V$ and it only has homology in dimensions at most $(m-1)(\dim(V)-1)$. Since $m\le k$, the product on the right hand side of the above formula is in fact finite (more precisely, it is non-zero only for $i=0, \dim(V)-1,2(\dim(V)-1),\ldots,(k-1)(\dim(V)-1)$). Therefore, we may think of the product as a direct sum, and so tensoring with $\R$ commutes with product in the displayed formulas below. 
By \refP{T_kCommutesWithReals}, we know that tensoring with $\R$ commutes, in our case, with holim, and so we obtain the weak equivalence
$$\left(T_k\Sing (\Ebar(M,V))\right)\otimes \R \simeq \left(\prod_{i=0}^\infty
\underset{(\phi,m)\in \tOsk(M)}{\holim}
\HH_i(\Balls(m,V);\Q)\right)\otimes \R$$
It is well-known (and is easy to prove using calculus of functors) that spaces such as $\Ebar(M,V)$ are homologically of finite type, therefore all chain complexes involved are homologically of finite type. Two rational chain complexes of homologically finite type that are quasi-isomorphic after tensoring with $\R$ are, necessarily, quasi-isomorphic over $\Q$. Therefore, we have a weak equivalence in $\CC_{\Q}$.
$$T_k\Sing (\Ebar(M,V)) \simeq \prod_{i=0}^\infty
\underset{(\phi,m)\in \tOsk(M)}{\holim}
\HH_i(\Balls(m,V);\Q)$$
The desired result follows by using, once again, \refP{compareEmbpiBLambda} and the equivalence between $\CC_{\Q}$ and rational spectra.

\end{proof}


\section{Formality and the splitting of the orthogonal tower}\label{S:FormalityOrthogonal}


In this section we show that \refT{LimPostnikovSplits}, which is about the splitting of a certain lim-Postnikov tower, can be reinterpreted as 
the splitting of the orthogonal tower of $\HQ\wedge\Ebar(M,V)_+$. Thus in this
section we mainly focus on the functoriality of $\HQ\wedge\Ebar(M,V)_+$ in $V$ and, accordingly, terms like ``Taylor polynomials'', ``derivatives'', etc. are always used in the context of orthogonal calculus\footnote{We are committing a slight abuse of notation here, because the definition of $\Ebar(M,V)$ depends on choosing a fixed embedding $M\hookrightarrow W$, and therefore $\Ebar(M,V)$ is only defined for vector spaces containing $W$. One way around this problem would be to work with the functor $V\mapsto \Ebar(M,W\oplus V)$. To avoid introducing ever messier notation, we chose to ignore this issue, as it does not affect our arguments in the slightest.}

As we have seen, embedding calculus tells us, roughly speaking, that  
$\Sigma^{\infty}\overline{\E}(M,V)_+$ can be written as a homotopy
inverse limit of spectra of the form $\Sigma^\infty\Config(k,V)_+$
where $\Config(k,V):=\Emb(\{1,\ldots, k\},V)$ is the space of configurations of $k$ points
in $V$. A good place to start is therefore to understand the
orthogonal Taylor tower of $V\mapsto\Sigma^\infty\Config(k,V)_+$.
The only thing that we will need in this section is the following
simple fact (we will only use a rationalized version of it, but it
is true integrally).
\begin{prop}\label{P:spherical}
The functor $V\mapsto\Sigma^\infty\Config(k,V)_+$ is polynomial of
degree $k-1$. For $0\le i\le k-1$, the $i$-th
layer in the orthogonal tower of this functor\footnote{Note that we are not speaking of the derivative of this functor}, $D_i\Sigma^\infty \Config(k,V)_+$,
is equivalent to a wedge of spheres of dimension $i(\dim(V)-1)$.
\end{prop}
This proposition is an immediate consequence of
\refP{LayersOfConfigurations} below, and its rational version is restated more precisely as \refC{ConfigurationsTower}.

We now digress to do a detailed calculation of the derivatives of $\Sigma^\infty \Config(k,V)_+$. First, we need some definitions.
\begin{definition}
Let $S$ be a finite set. A {\em partition} $\Lambda$ of $S$ is an
equivalence relation on $S$. Let $P(S)$ be the poset of all
partitions of $S$, ordered by refinement (the finer the bigger).
We say that a partition $\Lambda$ is {\em irreducible} if each
component of $\Lambda$ has at least $2$ elements.
\end{definition}

The poset $P(S)$ has both an initial and a final object. Therefore, its geometric realization, 
which we denote $|P(S)|$, is a contractible
simplicial complex. One standard way to construct a non-trivial homotopy type out of $P(S)$ is to consider the subposet $P_0(S)$, obtained from $P(S)$ by removing both the initial and the final object. We are going to construct a variation of $|P_0(S)|$ as follows. First, consider
the following sub-complex of $|P(S)|$, which we denote $\partial|P(S)|$. $\partial|P(S)|$ is built 
of those simplices of the nerve of $P(S)$ that do not contain the
morphism from the initial object to the final object as a $1$-dimensional face. It is an easy exercise to show that $\partial|P(S)|$ is homeomorphic to the unreduced suspension of $|P_0(S)|$.
Let $T_S$ be the quotient space $|P(S)|/\partial|P(S)|$. Since $|P(S)|$ is contractible, it follows that $T_S$ is equivalent to the suspension of $\partial|P(S)|$, and therefore to the double suspension of $|P_0(S)|$.

If $S=\{1,\ldots,n\}$, 
we denote $P(S)$ by $P(n)$ and $T_S$ by $T_n$.
There is a well-known equivalence (see, e.g., \cite[4.109]{OT}, where the analogous statement is proved for $|P_0(S)|$, which is called in [op. cit] the Folkman complex of the braid arrangement),
$$
T_n\simeq \bigvee_{(n-1)!}S^{n-1}.
$$

Now let $\Lambda$ be a partition of $S=\{1,\ldots,n\}$, and let $P(\Lambda)$ be the
poset of all refinements of $\Lambda$. Again, $P(\Lambda)$ is a poset with both an initial and a final object. Define $\partial|P(\Lambda)|$ in the same way as before, and let $$T_\Lambda=|P(\Lambda)|/\partial|P(\Lambda)|.$$ It
is not hard to see that if $\Lambda$ is a partition with
components $(\lambda_1,\ldots,\lambda_j)$ then there is an
isomorphism of posets
$$P(\Lambda)\cong P(\lambda_1)\times\cdots\times P(\lambda_j)$$
and therefore  a homeomorphism
$$T_\Lambda\cong T_{\lambda_1}\wedge \cdots\wedge T_{\lambda_j}.$$
In particular, $T_\Lambda$ is equivalent to a wedge of spheres of
dimension $n-j$. We call this number the {\em excess} of $\Lambda$
and denote it by $e(\Lambda)$.

\begin{prop}\label{P:LayersOfConfigurations}
For $i>0$, the $i$-th layer of $\Sigma^\infty \Config(k,V)_+$ is
equivalent to
$$D_i\Sigma^\infty\Config(k,V)_+\simeq\bigvee_{\{\Lambda\in P(k)\mid e(\Lambda)=i\}} \map_*\left(T_\Lambda,\Sigma^\infty S^{iV}\right)$$
where the wedge sum is over the set of partitions of $k$ of excess
$i$.
\end{prop}
\begin{proof}
Denote the fat diagonal of $kV$ by
$\Delta^kV:=\{(v_1,\cdots,v_k)\in kV\colon v_i=v_j \textrm{ for
some }i\not=j\}$. The smashed-fat-diagonal of $S^{kV}$ is
$$
\Delta^kS^V:=\{x_1\wedge\cdots\wedge x_k\in
\wedge_{i=1}^kS^V=S^{kV}: x_i=x_j \textrm{ for some }i\not=j\}.
$$
Thus
$$
\Config(k,V)=kV\setminus\Delta^kV=((kV)\cup\{\infty\})\setminus((\Delta^kV)\cup\{\infty\})
=S^{kV}\setminus\Delta^kS^V.$$

Recall that for a subpolyhedron in a sphere, $j\colon
K\hookrightarrow S^n$, Spanier-Whitehead duality gives a weak
equivalence of spectra
$$\Sigma^\infty(S^n\setminus K)_+\simeq\map_*(S^n/K,\Sigma^\infty
S^n)$$ which is natural with respect to inclusions $L\subset K$
and commutes with suspensions. In our case Spanier-Whitehead
duality gives  an equivalence
$$\Sigma^\infty \Config(k,V)_+ \simeq \map_*(S^{kV}/\Delta^kS^V,\Sigma^{\infty} S^{kV})$$
which is natural with respect to linear isometric injections. The
right hand side is equivalent to the homotopy fiber of the map
$$\map_*(S^{kV},\Sigma^\infty S^{kV})\longrightarrow \map_*(\Delta^kS^V,\Sigma^\infty S^{kV})$$
Since $\map_*(S^{kV},\Sigma^\infty S^{kV})\simeq \Sigma^\infty
S^0$ is a constant functor, it has no layers of degree greater than zero.
Therefore, for $i>0$,
$$D_i\Sigma^\infty \Config(k,V)_+ \simeq \Omega\left(D_i \map_*(\Delta^kS^V,\Sigma^\infty S^{kV})\right)$$
It is not hard to see (see \cite{SchwartzGenus}, Lemma 2.2 for a
proof) that $\Delta^kS^V$ can be ``filtered'' by excess. More precisely, there is a sequence of
spaces $$*=\Delta_0^kS^V\longrightarrow\Delta_1^kS^V\longrightarrow
\Delta_2^kS^V\longrightarrow\cdots\longrightarrow\Delta_{k-1}^kS^V=\Delta^k S^V$$ such
that the homotopy cofiber of the map $\Delta_{i-1}^kS^V\to
\Delta_{i}^kS^V$ is equivalent to
$$\bigvee_{\{\Lambda\in P(k)\mid e(\Lambda)=i\}} K_\Lambda \wedge S^{(k-i)V}$$
where $K_\Lambda$ is a de-suspension of $T_\Lambda$. It follows
that $\map_*(\Delta^kS^V,\Sigma^\infty S^{kV})$ can be decomposed
into a finite tower of fibrations
$$\map_*(\Delta^kS^V,\Sigma^\infty S^{kV})=X_{k-1}\longrightarrow X_{k-2}\longrightarrow\cdots\longrightarrow X_1$$
where the homotopy fiber of the map $X_i\to X_{i-1}$ is equivalent
to
$$\prod_{\{\Lambda\in P(k)\mid e(\Lambda)=i\}}\map_*(K_{\Lambda},\Sigma^\infty S^{iV})$$
Since this is obviously a homogeneous functor of degree $i$, it
follows that $X_i$ is the $i$-th Taylor polynomial of
$\map_*(\Delta^kS^V,\Sigma^\infty S^{kV})$. The proposition 
follows.
\end{proof}

Rationalizing, we obtain the following corollary.
\begin{cor}\label{C:ConfigurationsTower}
Each layer in the orthogonal tower of the functor
$V\mapsto\HQ\wedge\Config(k,V)_+$ is an Eilenberg-Maclane
spectrum.  More precisely,
$$ D_i(\HQ\wedge\Config(k,V)_+)\simeq
  \begin{cases}
    ||\HH_{i(\dim(V)-1)}(\Config(k,V))||  & \text{if }i\leq k-1; \\
    *, & \text{otherwise}
  \end{cases}
$$
where $||\HH_{i(\dim(V)-1)}(\Config(k,V))|| $ is the Eilenberg-Mac Lane spectrum
that has the $i(\dim(V)-1)$-th rational homology of $\Config(k,V)$ in degree
$i(\dim(V)-1)$.

Therefore, this orthogonal tower coincides, up to indexing,
with the Postnikov tower, i.e.
$$P_n(\HQ\wedge\Config(k,V)_+)\simeq\Po_{d(n)}(\HQ\wedge\Config(k,V)_+),$$
where $d(n)$ is any number satisfying $n(\dim V-1)\leq d(n)<(n+1)(\dim V-1)$.
\end{cor}
\begin{proof}
The computation of the layers is an immediate application of the
previous proposition. Set $X=\HQ\wedge\Config(k,V)_+$ and consider
the following commutative square
$$\xymatrix{X\ar[r]\ar[d]&\Po_d(X)\ar[d]\\P_n(X)\ar[r]&\Po_dP_n(X).}
$$
A study of the homotopy groups of the layers shows that the bottom
and the right maps are weak equivalences when $d$ is the
prescribed range.
\end{proof}

Let $\widehat{\Sigma^\infty\Config(k, V)_+}$ be the total homotopy fiber of the 
$k$-dimensional cubical diagram which sends a subset $S\subset\{1,\ldots, k\}$
to $\Sigma^\infty\Config(S, V)_+$ (where $\Config(S, V)=\Emb(S, V)$), and where the maps are the obvious restriction maps.
$\widehat{\Sigma^\infty\Config(k, V)_+}$ is a functor of $V$, and so we may ask about
the homogeneous layers of this functor. We have the following variation of Proposition~\ref{P:LayersOfConfigurations}.
\begin{prop}\label{P:LayersOfReducedConfigurations}
Let $P^{\mbox{irr}}(k)$ be the set of {\em irreducible} partitions of $k$ (i.e., partitions without singletons).
For $i\ge0$, the $i$-th layer of $\widehat{\Sigma^\infty \Config(k,V)_+}$ is
equivalent to
$$D_i\widehat{\Sigma^\infty\Config(k,V)_+}\simeq\bigvee_{\{\Lambda\in P^{\mbox{irr}}(k)\mid e(\Lambda)=i\}} \map_*\left(T_\Lambda,\Sigma^\infty S^{iV}\right)$$
where the wedge sum is over the set of irreducible partitions of $k$ of excess
$i$.
\end{prop}
\begin{proof}
It is clear by inspection that $$D_0\widehat{\Sigma^\infty\Config(k,V)_+}\simeq *$$ so we may assume that $i>0$. Proposition~\ref{P:LayersOfConfigurations} can be rephrased as saying that  for $i>0$, the $i$-th layer of $\Sigma^\infty \Config(S,V)_+$ is
equivalent to
$$\prod_{\{\Lambda\in P(S)\mid e(\Lambda)=i\}} \map_*\left(T_\Lambda,\Sigma^\infty S^{iV}\right)$$
where the product is over the set of partitions of $S$ of excess
$i$. Let $S_1\hookrightarrow S_2$ be an inclusion, and consider the corresponding projection of configuration spaces $\Config(S_2, V)\longrightarrow \Config(S_1, V)$. It is not hard to show, by inspecting the proof of Proposition~\ref{P:LayersOfConfigurations}, that the corresponding map 
of $i$-th layers
$$D_i\Sigma^\infty \Config(S_2, V)_+\longrightarrow D_i\Sigma^\infty \Config(S_1, V)_+$$
corresponds to the projection
$$\prod_{\{\Lambda\in P(S_2)\mid e(\Lambda)=i\}} \map_*\left(T_\Lambda,\Sigma^\infty S^{iV}\right)\longrightarrow \prod_{\{\Lambda\in P(S_1)\mid e(\Lambda)=i\}} \map_*\left(T_\Lambda,\Sigma^\infty S^{iV}\right)$$
associated with the obvious inclusion of posets $P(S_1)\hookrightarrow P(S_2)$. The proposition follows by elementary combinatorics.
\end{proof}
\begin{cor}
Suppose that $k=2l-1$ or $k=2l$. Then $D_i\widehat{\Sigma^\infty\Config(k,V)_+}\simeq *$ for
$i< l$.
\end{cor}
\begin{proof}
According to the preceding proposition, $D_i\widehat{\Sigma^\infty\Config(k,V)_+}$ is a wedge sum indexed by irreducible partitions of $k$ of excess $i$. It is clear from the definition of excess that the lowest possible excess of an irreducible partition of $k$ is attained by the irreducible partition of $k$ with the maximal number of components. It is easy to see that the irreducible partitions of $k$ with the maximal number of components are partitions of type $$\underbrace{2-\cdots-2}_{l-2}-3$$ if $k=2l-1$,  and partitions of type $$\underbrace{2-\cdots-2}_{l}$$ if $k=2l$. In both of these cases, the excess of the partition is $l$. It follows that the wedge sum in the statement of Proposition~\ref{P:LayersOfReducedConfigurations} is empty for $i<l$, and so the left hand side is contractible in this case.
\end{proof}
\begin{cor}\label{C:HighlyConnected}
Suppose that $k=2l-1$ or $k=2l$. Then $\widehat{\Sigma^\infty\Config(k,V)_+}$ is $l(\dim(V)-1)-1$-connected.
\end{cor}
\begin{proof}
 By the preceding corollary, the smallest $i$ for which $D_i\widehat{\Sigma^\infty\Config(k,V)_+}$ is non-trivial is $l$. By Proposition~\ref{P:LayersOfReducedConfigurations}, this layer is equivalent to a (stable) wedge of spheres of dimension $l(\dim(V)-1)$, and so it is $l(\dim(V)-1)-1$-connected. Clearly, higher layers are more highly connected. Since the taylor tower of $\Sigma^\infty\Config(k, V)_+$ converges, the statement follows.
\end{proof}
We will also need the following proposition.

\begin{prop}\label{P:OrthAndEmbCommute}
For every $n$ there exists a large enough $k$ such that the natural map 
$$P_n\HQ\wedge\Ebar(M,V)_+ \stackrel{\simeq}{\longrightarrow} \underset{U\in\mathcal{O}_k(M)}{\holim} P_n\HQ\wedge \Ebar(U,V)_+$$
is an equivalence.  The same holds if $P_n$ is replaced by $D_n$.
\end{prop}
\begin{proof}
 We will only prove the $P_n$
version. The target of the map is $T_kP_n\HQ\wedge \Ebar(M,V)_+$. Applying \refL{PnAndTkCommute} to the functor $E(M,V)=\HQ\wedge \Ebar(M,V)_+$, it is enough to prove that for a large enough $k$
the map $P_n\HQ\wedge \Ebar(M,V)_+\to P_nT_k\HQ\wedge
\Ebar(M,V)_+$ is an equivalence. Consider again the formula for the
$k$-th layer in the embedding tower
$$\Gamma_c\left({M\choose k}, \widehat{\HQ\wedge \Config(k,V)_+}\right).$$
It follows from Corollary~\ref{C:HighlyConnected} that the spectrum $ \widehat{\HQ\wedge \Config(k,V)_+}$ is roughly $\frac{k\dim(V)}{2}$-connected.  
It follows that $\HQ\wedge \Ebar(M,V)_+$ and $T_k\HQ\wedge\Ebar(M,V)_+$ agree to order roughly $\frac{k}{2}$ (in the sense defined in \refS{ReviewCalculi}). It follows, by \refP{ConnectivityEstimates}, that the map $\HQ\wedge \Ebar(M,V)_+\to T_k\HQ\wedge\Ebar(M,V)_+$ induces an equivalence on $P_n$, for roughly $n\le\frac{k}{2}$. \end{proof}
We are now ready to state and prove our main theorem

\begin{thm}\label{T:TowerSplits}
Under the assumptions of \refT{LimPostnikovSplits}, the
orthogonal tower of the functor $\HQ\wedge\overline{\E}(M,V)_+$
splits.
In other words, there is an equivalence
$$P_n\HQ\wedge \Ebar(M,V)_+\simeq\prod_{i=0}^n D_i\HQ\wedge \Ebar(M,V)_+.$$
\end{thm}

\begin{proof}
By \refL{PnAndTkCommute} and \refP{OrthAndEmbCommute}, and using
the model for $T_k\HQ\wedge \Ebar(M,V)_+$ given in \refT{LimPostnikovSplits}, it
is enough to show that
$$P_n\left(\underset{(m,\phi)\in \tOsk(M)}{\holim} \HQ\wedge \Balls(m,V)_+\right)\simeq \prod_{i=0}^n \underset{(m,\phi)\in \tOsk(M)}{\holim} D_i(\HQ\wedge \Balls(m,V)_+).$$
By \refC{ConfigurationsTower}  the Taylor tower of $\HQ\wedge
\Balls(m,V)_+$ coincides, up to regrading, with the Postnikov
tower. By the proof of \refT{LimPostnikovSplits}, the homotopy limit $\holim
\HQ\wedge\Balls(m,V)_+$ splits as a product of the homotopy limits
of the layers in the Postnikov towers. Since diagrams of layers in
the Postnikov towers and diagrams of layers in the orthogonal
towers are diagrams of Eilenberg-MacLane spectra that are equivalent on homotopy groups, they are
equivalent diagrams (as per Remark~\ref{R:UniquenessEMdiagrams}). It follows that $\holim \HQ\wedge\Balls(n,V)_+$
splits as a product of the homotopy limits of the layers in the
orthogonal towers.
\end{proof}
It is easy to see that the splitting is natural with respect to
inclusions  of submanifolds of $M$, but notice that we do not claim that the
splitting is natural in $V$.


\section{The layers of the orthogonal tower} \label{S:DescriptionOfLayers}


In this section we explicitly describe the layers (in the
sense of orthogonal calculus) of the Taylor tower of  $\HQ\wedge
\Ebar(M,V)$ as the twisted cohomology of certain spaces of
partitions attached to $M$. We will try to give a ``plausibility
argument'' for our formulas, but a detailed proof will appear
in \cite{Arone}.

We encountered partition posets in Section \ref{S:FormalityOrthogonal}. Here, however, we need to consider a different category of partitions. If
$\Lambda$ is a partition of $S$, we call $S$ the support of
$\Lambda$. When we need to emphasize that $S$ is the support of $\Lambda$,
we use the notation $S(\Lambda)$. Also, we denote by $c(\Lambda)$ the set of components of
$\Lambda$. Then $\Lambda$ can be represented by a surjection
$S(\Lambda) \twoheadrightarrow c(\Lambda)$. Let $C_\Lambda$ be the mapping cylinder of this surjection. Then $S(\Lambda)\subset C_\Lambda$. In the previous
section we defined the excess of $\Lambda$ to be $e(\Lambda):=
|S(\Lambda)|- |c(\Lambda)|$. It is easy to see that
$$e(\Lambda)=\mbox{rank}(\HH_1(C_\Lambda,S(\Lambda)).$$

Let $\Lambda_1, \Lambda_2$ be partitions of $S_1, S_2$
respectively. A ``pre-morphism'' $\alpha:\Lambda_1\to \Lambda_2$
is defined to be a surjection (which we denote with the same
letter) $\alpha:S_1 \twoheadrightarrow S_2$ such that $\Lambda_2$ is
the equivalence relation generated by $\alpha(\Lambda_1)$. It is
easy to see that such a morphisms induces a map of pairs $(C_{\Lambda_1},S(\Lambda_1))\to (C_{\Lambda_2},S(\Lambda_2))$, and therefore a homomorphism
$$\alpha_*:\HH_1(C_{\Lambda_1},S(\Lambda_1))\longrightarrow
\HH_1(C_{\Lambda_2},S(\Lambda_2)).$$ We say that $\alpha$ is a
{\em morphism} if $\alpha_*$ is an isomorphism. In particular,
there can only be a morphism between partitions of
equal excess. Roughly speaking, morphisms are allowed to fuse components
together, but are not allowed to bring together two elements in
the same component.

For $k\ge 1$, let ${\mathcal E}_k$ be the category of irreducible
partitions (recall that $\Lambda$ is irreducible if none of the
components of $\Lambda$ is a singleton) of excess $k$, with
morphisms as defined above. Notice that if $\Lambda$ is irreducible
of excess $k$ then the size of the support of $\Lambda$ must be
between $k+1$ and $2k$.

Next we define two functors on ${\mathcal E}_k$ -- one
covariant and one contravariant. Recall from the previous section
that $P(\Lambda)$ is the poset of refinements of $\Lambda$. A
morphism $\Lambda\to\Lambda'$ induces a map of posets
$P(\Lambda)\to P(\Lambda')$.  It is not difficult to see that this
map takes boundary into boundary, and therefore it induces a map
$T_{\Lambda}\to T_{\Lambda'}$. This construction
 gives rise to a functor ${\mathcal
E}_k\to\Top$, given on objects by
$$ \Lambda \longmapsto T_\Lambda.$$

In fact, to conform with the classification of homogeneous
functors in orthogonal calculus, we would like to induce up
$T_\Lambda$ to make a space with an action of the orthogonal group
$\OO(k)$. Let
$$\widetilde T_\Lambda:=\iso(\mathbb{R}^k,\HH_1(T(\Lambda),S(\Lambda);\mathbb{R}))_+\wedge T_\Lambda$$ where $\iso(V,W)$ is the space of linear isomorphisms from $V$ to $W$ (thus $\iso(V,W)$ is abstractly homeomorphic to the general linear group if $V$ and $W$ are isomorphic, and is  empty otherwise). 
In this way we get a covariant functor from ${\mathcal E}_k$ to spaces with
an action of $\OO(k)$.

Now we construct another functor on ${\mathcal E}_k$, this time a contravariant
one. To begin with, there is an obvious contravariant functor, defined
on objects by
$$\Lambda\mapsto M^{S(\Lambda)}.$$
Let $f\colon S(\Lambda)\longrightarrow M$ be an element of $M^{S(\Lambda)}$. 
The image of $f$ is a finite subset of $M$, and $f(\Lambda)$ is a partition of 
$f(S(\Lambda))$. Clearly, $f$ defines a pre-morphism $\Lambda\longrightarrow f(\Lambda)$.
We say that $f$ is a {\em good} element of $M^{S(\Lambda)}$ if the pre-morphism $\Lambda\longrightarrow f(\Lambda)$ is in fact a morphism. Otherwise we say that $f$ is a bad element of $M^{S(\Lambda)}$.
\begin{eg}
Let $\Lambda$ be the partition $(1,2)(3,4)$. Let $f$ be a map $f\colon \{1,2,3,4\} \longrightarrow M$.
If $f$ is injective, then $f$ is good. If $f(2)=f(3)$, but otherwise $f$ is injective, then $f$ is good. If $f(1)=f(2)$ then $f$ is bad. In general, if $f$ is not injective on some component of $\Lambda$, then $f$ is bad, but the converse is not true. In our example, if $f(1)=f(3)$ and $f(2)=f(4)$ then $f$ is bad, even though it may be injective on each component.
\end{eg}
Let $\Delta^\Lambda(M)$ be the suspace of $M^{S(\Lambda)}$ consisting of all the bad elements. 
For example, if $\Lambda$
is the partition with one component then $\Delta^\Lambda(M)$ is
the usual fat diagonal. It is not hard to see that the contravariant functor
$\Lambda\mapsto M^{S(\Lambda)}$ passes to a
contravariant functor from ${\mathcal E}_k$ to spaces given on objects by
$$\Lambda\mapsto M^{S(\Lambda)}/\Delta^\Lambda(M).$$
Let $M^{[\Lambda]}:=M^{S(\Lambda)}/\Delta^\Lambda(M)$.
Given a covariant functor and a contravariant functor on ${\mathcal E}_k$, we may consider the ``tensor product'' (homotopy coend)
$$\widetilde T_\Lambda\otimes_{{\mathcal E}_k} M^{[\Lambda]},$$
which is a space with an action of $\OO(k)$.
\begin{thm}\label{T:Description}
The $i$-th layer of the orthogonal calculus tower of
$\Sigma^\infty \Ebar(M,V)_+$ is equivalent to
$$\map_*\left(\widetilde T_\Lambda\otimes_{{\mathcal E}_i} M^{[\Lambda]},\Sigma^\infty S^{Vi}\right)^{\OO(i)}.$$
\end{thm}
\begin{proof}[Idea of proof]
Embedding calculus suggests that it is almost enough to prove the
theorem in the case of $M$ homeomorphic to a finite disjoint union
of balls. In this case $\Ebar(M,V)$ is equivalent to the
configuration space $\Config(k,V)$. It is not hard to show that then the formula in the statement of the current theorem is
equivalent to the formula given by \refP{LayersOfConfigurations}.
The current theorem restates the formula of
\refP{LayersOfConfigurations} in a way that is well-defined and natural for
all manifolds $M$.
\end{proof}
It follows that the $k$-th layer of $\HQ\wedge \Ebar(M,V)_+$ is given by the same formula as in the theorem, with $\Sigma^\infty$ replaced with $\HQ\wedge$.
\begin{cor}\label{C:RHT-invariance-layers}
Suppose that $M_1$ and $M_2$ are related by a zig-zag of map inducing an isomorphism in homology. Then for each $n$, the $n$-th layers of the orthogonal towers of the two functors 
$$V\mapsto
\Sigma^\infty\Ebar(M_i,V)_+, \quad i=1,2$$ are homotopy equivalent. Similarly, if the maps in the aforementioned zig-zag induce a isomorphisms in rational homology then the layers of the orthogonal towers of $V\mapsto
\HQ\wedge\Ebar(M_i,V)_+$ are equivalent.
\end{cor} 
\begin{proof}
It is not hard to show that $\widetilde T_\Lambda\otimes_{{\mathcal E}_k} M^{[\Lambda]}$ is equivalent to a finite CW complex with a free action (in the pointed sense) of $\OO(k)$. Since the action is free, the fixed points construction in the formula for the layers in the orthogonal tower can be replaced with the homotopy fixed points construction. Thus, the $k$-th layer in the orthogonal tower of $\Sigma^\infty \Ebar(M,V)_+$ 
is equivalent to  
$$\map_*\left(\widetilde T_\Lambda\otimes_{{\mathcal E}_k} M^{[\Lambda]},\Sigma^\infty S^{Vk}\right)^{h\OO(k)}.$$
It is easy to see that this is a functor that takes homology equivalences in $M$ to homotopy equivalences. For the rational case, notice
that $$\map_*\left(\widetilde T_\Lambda\otimes_{{\mathcal E}_k} M^{[\Lambda]},\HQ\wedge S^{Vk}\right)^{h\OO(k)}$$ is a functor of $M$ that takes rational homology equivalences to homotopy equivalences.
\end{proof}

 Some remarks are perhaps in order.

\begin{rem}  It may be helpful to note that the space
$\widetilde T_\Lambda\otimes_{{\mathcal E}_k} M^{[\Lambda]}$ can
be filtered by the size of support of $\Lambda$ (that is, by the
number of points in $M$ involved). This leads to a decomposition
of the $k$-th layer in the orthogonal tower of
$\Sigma^\infty\Ebar(M,V)$ as a finite tower of fibrations, with
$k$ terms, indexed $k+1\le i\le 2k$, corresponding to the number
of points in $M$. This is the embedding tower of the $k$-th layer
of the orthogonal tower. For example, the second layer of the
orthogonal tower fits into the following diagram, where $\Delta^{2,2}M$ is the singular
set of the action of $\Sigma_2\wr \Sigma_2$ on $M^4$, the left
row is a fibration sequence, and the square is a homotopy
pullback.

$$
\xymatrix{ \map_*(\frac{M^4}{\Delta^4M}\wedge T_2\wedge
T_2,\Sigma^\infty S^{2V})^{\Sigma_2\wr\Sigma_2}
\ar[d]  &    \\
 D_2\Sigma^\infty\Ebar(M,V)
 \ar[r] \ar[d]  &
\map_*(\frac{M^4}{\Delta^{2,2}M}\wedge T_2^{\wedge
2},\Sigma^\infty S^{2V})^{\Sigma_2\wr\Sigma_2}
\ar[d]  \\
\map_*(\frac{M^3}{\Delta^3M}\wedge T_3,\Sigma^\infty
S^{2V})^{\Sigma_3 } \ar[r]  & \map_*(\frac{M^3}{\Delta^3M}\wedge
T_2^{\wedge 2},\Sigma^\infty S^{2V})^{\Sigma_2} }
$$
\end{rem}


\begin{rem}
 To relate this to something ``classical'', note that the top layer of the embedding
 tower of the $k$-th layer of the orthogonal tower is
$$\map_*(M^{2k}/\Delta^{2k}M\wedge T_2^{\wedge k},\Sigma^\infty S^{kV})^{\Sigma_2\wr\Sigma_k}.$$
This is the space of ``chord diagrams'' on $M$, familiar from knot
theory. In fact, in the case of $M$ being a circle (or an interval, in which case one considers embeddings fixed near the boundary), it is known from \cite{LV} that the Vassiliev homology spectral sequence, which also converges to the space of knots, collapses at $E^1$.  Thus 
the orthogonal tower spectral sequence for $\HQ\wedge \Ebar(M,V)$ must coincide
with Vassiliev's.  It is not hard to verify directly that the two $E^1$ terms are isomorphic (up to regrading).
\end{rem}

Finally, we deduce the rational homology invariance of
$\Ebar(M,V)$ from our main theorem and Corollary~\ref{C:RHT-invariance-layers}.
\begin{thm} \label{T:RationalHomologyInvariant}
Let $M$ and $M'$ be two manifolds such that there is
a zig-zag of maps, each inducing an isomorphism in rational homology, connecting  $M$ and $M'$.
If $$\dim V > 2\max(\ed( M),\ed (M')),$$
then $\Ebar(M,V)$ and $\Ebar(M',V)$ have the same rational
homology groups.
\end{thm}

\end{document}